\documentclass[11pt, a4paper]{article}
\usepackage{amsmath}
\usepackage{wasysym}
\usepackage{latexsym}
\usepackage{amsfonts}
\usepackage{mathrsfs}
\usepackage{amssymb}
\usepackage{ifsym}
\usepackage{dsfont}
\usepackage[all]{xy}
\usepackage{authblk}

\newtheorem{Definitions1}{Definition}[section]

\newtheorem{Theorems1}{Theorem}[section]

\newtheorem{Coroll1}[Theorems1]{Corollary}
\newtheorem{Lemma1}[Theorems1]{Lemma}

\newtheorem{Remark1}[Theorems1]{Remark}

\newtheorem{Quest1}{Question}[section]

\newenvironment{proof}[1][Proof]{\begin{trivlist}
\item[\hskip \labelsep {\bfseries #1}]}{\end{trivlist}}

\newcommand{\qed}{\nobreak \ifvmode \relax \else
      \ifdim\lastskip<1.5em \hskip-\lastskip
      \hskip1.5em plus0em minus0.5em \fi \nobreak
      \vrule height0.75em width0.5em depth0.25em\fi}

\begin{document}
\title{On the relative strengths of fragments of collection} 
\author[1]{Zachiri McKenzie}
\affil[1]{University of Michigan- Shanghai Jiao Tong University Joint Institute,
Shanghai Jiao Tong University, Shanghai, P. R. China\newline
\texttt{zach.mckenzie@gmail.com}}
\maketitle

\begin{abstract}
Let $\mathbf{M}$ be the basic set theory that consists of the axioms of extensionality, emptyset, pair, union, powerset, infinity, transitive containment, $\Delta_0$-separation and set foundation. This paper studies the relative strength of set theories obtained by adding fragments of the set-theoretic collection scheme to $\mathbf{M}$. We focus on two common parameterisations of the collection: {\it $\Pi_n$-collection}, which is the usual collection scheme restricted to $\Pi_n$-formulae, and {\it strong $\Pi_n$-collection}, which is equivalent to $\Pi_n$-collection plus $\Sigma_{n+1}$-separation. The main result of this paper shows that for all $n \geq 1$,
\begin{enumerate}
\item $\mathbf{M}+\Pi_{n+1}\textrm{-collection}+\Sigma_{n+2}\textrm{-induction on } \omega$ proves that there exists a transitive model of Zermelo Set Theory plus $\Pi_{n}$-collection, 
\item the theory $\mathbf{M}+\Pi_{n+1}\textrm{-collection}$ is $\Pi_{n+3}$-conservative over the theory $\mathbf{M}+\textrm{strong }\Pi_n \textrm{-collection}$.
\end{enumerate}
It is also shown that (2) holds for $n=0$ when the Axiom of Choice is included in the base theory. The final section indicates how the proofs of (1) and (2) can be modified to obtain analogues of these results for theories obtained by adding fragments of collection to a base theory (Kripke-Platek Set Theory with Infinity plus $V=L$) that does not include the powerset axiom. 
\end{abstract}

\section[Introduction]{Introduction}

In \cite{mat01} Mathias systematically studies and compares a variety of subsystems of $\mathrm{ZFC}$. One of the weakest systems studied in \cite{mat01} is the set theory $\mathbf{M}$ axiomatised by: extensionality, emptyset, pair, union, powerset, infinity, transitive containment, $\Delta_0$-separation and set foundation. This paper will expand upon some of the initial comparisons of extensions of $\mathbf{M}$ achieved in \cite{mat01} by studying the strengths of extensions of $\mathbf{M}$ obtained by adding fragments of the set-theoretic collection scheme. The fragments of the collection scheme considered in this paper will be obtained by restricting the following alternative versions of the collection scheme to the Takahashi class $\Delta_0^\mathcal{P}$ and the L\'{e}vy $\Pi_n$ classes:
\begin{itemize}
\item[](Collection) For all formulae $\phi(x, y, \vec{z})$ in the language of set theory,
$$\forall \vec{z} \forall w((\forall x \in w) \exists y \phi(x, y, \vec{z}) \Rightarrow \exists C (\forall x \in w)(\exists y \in C) \phi(x, y, \vec{z})).$$ 
\item[](Strong Collection) For all formulae $\phi(x, y, \vec{z})$ in the language of set theory,
$$\forall \vec{z} \forall w \exists C (\forall x \in w)(\exists y \phi(x, y, \vec{z}) \Rightarrow (\exists y \in C) \phi(x, y, \vec{z})).$$
\end{itemize}
Both Collection and Strong Collection yield $\mathrm{ZF}$ when added to $\mathbf{M}$. In section \ref{Sec:Background} we note that, over $\mathbf{M}$, the restriction of the Strong Collection scheme to $\Pi_n$-formulae (strong $\Pi_n$-collection) is equivalent to the restriction of the Collection scheme to $\Pi_n$-formulae ($\Pi_n$-collection) plus separation for all $\Sigma_{n+1}$-formulae. This means that $\mathbf{M}$ plus $\Pi_{n+1}$-collection proves all instances of strong $\Pi_n$-collection.  

One of the many achievements of \cite{mat01} is showing that if $\mathbf{M}$ is consistent, then so is $\mathbf{M}$ plus the Axiom of Choice and strong $\Delta_0$-collection. In section \ref{Sec:BaseCaseSection} we investigate the strength of adding $\Delta_0^\mathcal{P}$-collection to four of the weak set theories studied in \cite{mat01}. We show that if $T$ is one of the theories $\mathbf{M}$, $\mathrm{Mac}$, $\mathbf{M}+\mathrm{H}$ or $\mathrm{MOST}$, then $T$ plus $\Delta_0^\mathcal{P}$-collection is $\Pi_2^\mathcal{P}$-conservative over $T$. As a consequence, we are able to extend the consistency results of \cite{mat01} by showing that if $\mathbf{M}$ is consistent, then so is $\mathbf{M}$ plus the Axiom of Choice plus $\Pi_1$-collection.

The results of \cite{mat01} also show that the theory obtained by adding strong $\Pi_1$-collection to $\mathbf{M}$ is strictly stronger than $\mathbf{M}$. More specifically, $\mathbf{M}$ plus strong $\Pi_1$-collection proves the consistency of Zermelo Set Theory plus $\Delta_0$-collection. This result and the main result of section \ref{Sec:BaseCaseSection} are generalised in section \ref{Sec:GeneralCollectionResults} to show: For all $n \geq 1$,
\begin{enumerate} 
\item $\mathbf{M}$ plus $\Pi_{n+1}$-collection and the scheme of induction on $\omega$ restricted to $\Sigma_{n+2}$-formulae proves that there exists a transitive model of Zermelo Set Theory plus $\Pi_n$-collection,
\item the theory $\mathbf{M}+\Pi_{n+1}\textrm{-collection}$ is $\Pi_{n+3}$-conservative over the theory\\ 
$\mathbf{M}+\textrm{strong }\Pi_n \textrm{-collection}$.
\end{enumerate}
These comparisons are achieved using techniques, developed by Pino and Ressayre in \cite{res87} (see also \cite{flw16}), for building models of fragments of the collection scheme from chains of partially elementary submodels of the universe indexed by an ordinal, or a cut of a nonstandard ordinal, of a model of set theory.  

Finally, in section \ref{Sec:ResultsForKP} we consider replacing the base theory $\mathbf{M}$ by a theory, Kripke-Platek Set Theory with the Axiom of Infinity ($\mathrm{KPI}$) plus $V=L$, that does not include the powerset axiom. We indicate how the arguments in section \ref{Sec:GeneralCollectionResults} can be adapted to obtain the following analogues of (1) and (2) above: For all $n \in \omega$,
\begin{enumerate}
\item $\mathrm{KPI}+V=L$ plus $\Pi_{n+1}$-collection and the scheme of induction on $\omega$ restricted to $\Sigma_{n+2}$-formulae proves that there exists a transitive model of the theory $\mathrm{KPI}+V=L$ plus strong $\Pi_n$-collection, and full class foundation,  
\item the theory $\mathrm{KPI}+V=L$ plus $\Pi_{n+1}$-collection is $\Pi_{n+3}$-conservative over the theory $\mathrm{KPI}+V=L$ plus strong $\Pi_{n}$-collection. 
\end{enumerate}

\section[Background]{Background} \label{Sec:Background}

Throughout this paper $\mathcal{L}$ will denote the language of set theory. Structures will usually be denoted using upper-case calligraphy roman letters ($\mathcal{M}, \mathcal{N}, \ldots$) and the corresponding plain font letter ($M, N, \ldots$) will be used to denote the underlying set of that structure. If $\mathcal{M}$ is a structure, then we will use $\mathcal{L}(\mathcal{M})$ to denote the language of $\mathcal{M}$. If $\mathcal{M}$ is an $\mathcal{L}^\prime$-structure where $\mathcal{L}^\prime \supseteq \mathcal{L}$ and $a \in M$ then we will use $a^*$ to denote the class $\{x \in M \mid \mathcal{M} \models (x \in a)\}$. As usual $\Delta_0 (= \Sigma_0= \Pi_0), \Sigma_1, \Pi_1, \ldots$ will be used to denote the L\'{e}vy classes of $\mathcal{L}$-formulae, and we use $\Pi_\infty$ to denote the union of all of these classes (i.e. $\Pi_\infty= \bigcup_{n \in \omega} \Sigma_n= \bigcup_{n \in \omega} \Pi_n$). For all $n \in \omega$, $\Delta_n$ is the class of all formulae that are provably equivalent to both a $\Sigma_n$ formula and a $\Pi_n$ formula. We will also have cause to consider the class $\Delta_0^\mathcal{P}$, which is the smallest class of $\mathcal{L}$-formulae that contains all atomic formulae, contains all compound formulae formed using the connectives of first-order logic, and is closed under quantification in the form $\mathcal{Q} x \in y$ and $\mathcal{Q} x \subseteq y$ where $x$ and $y$ are distinct variables, and $\mathcal{Q}$ is $\exists$ or $\forall$. The classes $\Sigma_1^\mathcal{P}, \Pi_1^\mathcal{P}, \Delta_1^\mathcal{P}, \ldots$ are defined inductively from the class $\Delta_0^\mathcal{P}$ in the same way that the classes $\Sigma_1, \Pi_1, \Delta_1, \ldots$ are defined from $\Delta_0$. If $\Gamma$ is a class of formulae and $T$ is a theory, then we write $\Gamma^T$ for the class of formulae that are provably equivalent in $T$ to a formula in $\Gamma$. If $\Gamma$ is a class of formulae, then we use $\mathbf{Bol}(\Gamma)$ to denote the smallest class of formulae that contains $\Gamma$, and contains all compound formulae formed using the connectives of first-order logic. Note that for all $n \in \omega$, $\mathbf{Bol}(\Sigma_n)^\emptyset= \mathbf{Bol}(\Pi_n)^\emptyset$ and $\mathbf{Bol}(\Sigma_n^\mathcal{P})^\emptyset= \mathbf{Bol}(\Pi_n^\mathcal{P})^\emptyset$. If $\Gamma$ is a class of formulae, then we write $\neg \Gamma$ for the class of negations of formulae in $\Gamma$. So, for all $n \in \omega$, $(\neg \Sigma_n)^\emptyset= \Pi_n^\emptyset$, $(\neg \Pi_n)^\emptyset= \Sigma_n^\emptyset$, $(\neg \Sigma_n^\mathcal{P})^\emptyset=(\Pi_n^\mathcal{P})^\emptyset$, and $(\neg \Pi_n^\mathcal{P})^\emptyset= (\Sigma_n^\mathcal{P})^\emptyset$. Let $T$ be an $\mathcal{L}^\prime$-theory and let $S$ be $\mathcal{L}^{\prime\prime}$-theory where $\mathcal{L}^{\prime} \subseteq \mathcal{L}^{\prime\prime}$, and let $\Gamma$ be a class of $\mathcal{L}^\prime$-formulae. The theory $S$ is said to be {\it $\Gamma$-conservative} over $T$ if $S$ and $T$ prove the same $\Gamma$-sentences. 


Let $\mathcal{M}$ and $\mathcal{N}$ be $\mathcal{L}$-structures. If $\mathcal{M}$ is a substructure of $\mathcal{N}$ then we will write $\mathcal{M} \subseteq \mathcal{N}$. If $\Gamma$ is a class of $\mathcal{L}$-formulae then we will write $\mathcal{M} \prec_\Gamma \mathcal{N}$ if $\mathcal{M} \subseteq \mathcal{N}$ and for every $\vec{a} \in M$, $\vec{a}$ satisfies the same $\Gamma$-formulae in both $\mathcal{M}$ and $\mathcal{N}$. In the case that $\Gamma$ is $\Pi_\infty$ or  $\Sigma_n$ then we will abbreviate this notation by writing $\mathcal{M} \prec \mathcal{N}$ and $\mathcal{M} \prec_n \mathcal{N}$ respectively. If $\mathcal{M} \subseteq \mathcal{N}$ and for all $x \in M$ and $y \in N$,
$$\textrm{if } \mathcal{N} \models (y \in x) \textrm{ then } y \in M,$$
then we say that $\mathcal{N}$ is an \emph{end-extension} of $\mathcal{M}$ and write $\mathcal{M} \subseteq_e \mathcal{N}$. It is well-known that if $\mathcal{M} \subseteq_e \mathcal{N}$ then $\mathcal{M} \prec_0 \mathcal{N}$. The following is a slight generalisation of the notion of a powerset preserving end-extension that was first studied by Forster and Kaye in \cite{fk91}.

\begin{Definitions1}
Let $\mathcal{M}$ and $\mathcal{N}$ be $\mathcal{L}$-structures. We say that $\mathcal{N}$ is a powerset preserving end-extension of $\mathcal{M}$, and write $\mathcal{M} \subseteq_e^\mathcal{P} \mathcal{N}$ if 
\begin{itemize}
\item[(i)] $\mathcal{M} \subseteq_e \mathcal{N}$, 
\item[(ii)] for all $x \in N$ and for all $y \in M$, if $\mathcal{N} \models (x \subseteq y)$, then $x \in M$.
\end{itemize} 
\end{Definitions1} 

Just as end-extensions preserve $\Delta_0$ properties, powerset preserving end-extensions preserve $\Delta_0^\mathcal{P}$ properties. The following is a slight modification of a result proved in \cite{fk91}:

\begin{Lemma1}
Let $\mathcal{M}$ and $\mathcal{N}$ be $\mathcal{L}$-structures that satisfy extensionality. If $\mathcal{M} \subseteq_e^\mathcal{P} \mathcal{N}$, then $\mathcal{M} \prec_{\Delta_0^\mathcal{P}} \mathcal{N}$. \Square
\end{Lemma1} 

Let $\Gamma$ be a class of $\mathcal{L}$-formulae. The following define the restriction of some commonly encountered axiom and theorem schemes of $\mathrm{ZFC}$ to formulae in the class $\Gamma$:
\begin{itemize}
\item[]($\Gamma$-separation) For all $\phi(x, \vec{z}) \in \Gamma$,
$$\forall \vec{z} \forall w \exists y \forall x(x \in y \iff (x \in w) \land \phi(x, \vec{z})).$$  
\item[]($\Gamma$-collection) For all $\phi(x, y, \vec{z}) \in \Gamma$,
$$\forall \vec{z} \forall w((\forall x \in w) \exists y \phi(x, y, \vec{z}) \Rightarrow \exists C (\forall x \in w)(\exists y \in C) \phi(x, y, \vec{z})).$$ 
\item[](strong $\Gamma$-collection) For all $\phi(x, y, \vec{z}) \in \Gamma$,
$$\forall \vec{z} \forall w \exists C (\forall x \in w)(\exists y \phi(x, y, \vec{z}) \Rightarrow (\exists y \in C) \phi(x, y, \vec{z})).$$
\item[]($\Gamma$-foundation) For all $\phi(x, \vec{z}) \in \Gamma$,
$$\forall \vec{z}(\exists x \phi(x, \vec{z}) \Rightarrow \exists y(\phi(y, \vec{z}) \land (\forall x \in y) \neg \phi(x, \vec{z}))).$$
If $\Gamma= \{x \in z\}$ then we will refer to $\Gamma$-foundation as \emph{set foundation}.
\item[]($\Gamma$-induction on $\omega$) For all $\phi(x, \vec{z}) \in \Gamma$,
$$\forall \vec{z}(\phi(\emptyset, \vec{z}) \land (\forall n \in \omega)(\phi(n, \vec{z}) \Rightarrow \phi(n+1, \vec{z})) \Rightarrow (\forall n \in \omega) \phi(n, \vec{z})).$$ 
\end{itemize}
We will use $\bigcup x \subseteq x$ to abbreviate the $\Delta_0$-formula that says that $x$ is transitive ($(\forall y \in x)(\forall z \in y)(z \in x)$). We will also make reference to the following axioms:
\begin{itemize}
\item[](Axiom H)
$$\forall u \exists T \left(\bigcup T \subseteq T \land \forall z(\bigcup z \subseteq z \land |z| \leq |u| \Rightarrow z \subseteq T) \right).$$
\item[]($\mathrm{TCo}$) $$\forall x \exists y \left( \bigcup y \subseteq y \land x \subseteq y \right).$$ 
\end{itemize}

The following weak subsystems of $\mathrm{ZFC}$ are studied by Mathias in \cite{mat01}:
\begin{itemize}
\item $\mathbf{S}_1$ is the $\mathcal{L}$-theory with axioms: extensionality, emptyset, pair, union, set difference, and powerset.
\item $\mathbf{M}$ is obtained from $\mathbf{S}_1$ by adding $\mathrm{TCo}$, infinity, $\Delta_0$-separation, and set foundation.  
\item $\mathrm{Mac}$ is obtained from $\mathbf{M}$ by adding the axiom of choice.
\item $\mathbf{M}+\mathrm{H}$ is obtained from $\mathbf{M}$ by adding Axiom H.  
\item $\mathrm{KPI}$ is obtained from $\mathbf{M}$ by removing powerset, and adding $\Delta_0$-collection and $\Pi_1$-foundation.
\item $\mathrm{KP}^\mathcal{P}$ is obtained from $\mathbf{M}$ by adding $\Delta_0^\mathcal{P}$-collection and $\Pi_1^\mathcal{P}$-foundation.
\item $\mathrm{MOST}$ is obtained from $\mathrm{Mac}$ by adding $\Sigma_1$-separation and $\Delta_0$-collection.
\item $\mathrm{Z}$ is obtained from $\mathbf{M}$ by removing $\mathrm{TCo}$, and adding $\Pi_\infty$-separation.
\item $\mathrm{ZC}$ is obtained from $\mathrm{Z}$ by adding the axiom of choice.  
\end{itemize}
In addition to these theories, we will also use $\mathrm{MOST}^{-\mathrm{AC}}$ to refer to the theory obtained by removing the axiom of choice from $\mathrm{MOST}$, and $\mathbf{M}^-$ to refer to the theory obtained by removing the powerset axiom from $\mathbf{M}$. $\mathrm{ZF}$ and  $\mathrm{ZFC}$ are obtained by adding $\Pi_\infty$-collection (or, equivalently, strong $\Pi_\infty$-collection) to $\mathbf{M}$ and $\mathrm{Mac}$ respectively.   

We begin by collecting together some well-known relationships between fragments of induction, separation, collection, and strong collection over the weak base theory $\mathbf{M}^-$. 

\begin{Lemma1} \label{Th:BasicRelationships}
Let $\Gamma$ be a class of $\mathcal{L}$-formulae. Let $n \in \omega$.
\begin{enumerate}
\item $\mathbf{M}^- + \Gamma\textrm{-foundation} \vdash \neg \Gamma\textrm{-induction on } \omega$ 
\item $\mathbf{M}^- + \Gamma\textrm{-separation} \vdash \mathbf{Bol}(\Gamma)\textrm{-separation}$
\item $\mathbf{M}^- + \Gamma\textrm{-separation} \vdash \Gamma\textrm{-foundation}$
\item $\mathbf{M}^- + \textrm{[strong] } \Pi_n^{(\mathcal{P})}\textrm{-collection} \vdash \textrm{[strong] } \Sigma_{n+1}^{(\mathcal{P})}\textrm{-collection}$
\item $\mathbf{M}^- + \Pi_n\textrm{-collection} \vdash \Delta_{n+1}\textrm{-separation}$   
\end{enumerate}
\Square
\end{Lemma1}

Another well-known application of $\Pi_n$-collection is that, over $\mathbf{M}^-$, this scheme implies that the classes $\Sigma_{n+1}$ and $\Pi_{n+1}$ are essentially closed under bounded quantification.

\begin{Lemma1} \label{Th:Normalisation}
Let $\phi(x, \vec{z})$ be a $\Sigma_{n+1}$-formula, and let $\psi(x, \vec{z})$ be a $\Pi_{n+1}$-formula. The theory $\mathbf{M}^-+ \Pi_n\textrm{-collection}$ proves that $(\forall x \in y) \phi(x, \vec{z})$ is equivalent to a $\Sigma_{n+1}$-formula, and $(\exists x \in y) \psi(x, \vec{z})$ is equivalent to a $\Pi_{n+1}$-formula. \Square  
\end{Lemma1}  

We also observe that for all $n \in \omega$, strong $\Pi_n$-collection is equivalent, over $\mathbf{M}^-$, to $\Pi_n$-collection plus $\Sigma_{n+1}$-separation. The following lemma generalises one of the equivalences reported in \cite[Proposition 3.14]{mat01}. 

\begin{Lemma1} \label{Th:StrongCollectionVsCollectionSeparation}
For all $n \in \omega$,
\begin{enumerate}
\item $\mathbf{M}^- + \textrm{strong } \Pi_n\textrm{-collection} \vdash \Pi_n\textrm{-collection and } \Sigma_{n+1}\textrm{-separation}$ 
\item $\mathbf{M}^- + \Pi_n\textrm{-collection}+\Sigma_{n+1}\textrm{-separation} \vdash \textrm{strong } \Pi_n \textrm{-collection}$.
\end{enumerate}
\end{Lemma1}

\begin{proof}
We first prove (1). The fact that $\mathbf{M}^- + \textrm{strong } \Pi_n\textrm{-collection}$ proves the scheme of $\Pi_n$-collection is clear. We need to prove that $\mathbf{M}^- + \textrm{strong } \Pi_n\textrm{-collection}$ proves the scheme of $\Sigma_{n+1}$-separation. It immediately follows from Lemma \ref{Th:BasicRelationships} that $\mathbf{M}^- + \textrm{strong } \Pi_n\textrm{-collection}$ proves the scheme of strong $\Sigma_{n+1}$-collection and $\Pi_n$-separation. Work in the theory $\mathbf{M}^- + \textrm{strong } \Pi_n\textrm{-collection}$. Consider $\exists y \phi(y, x, \vec{z})$ where $\phi(y, x, \vec{z})$ is $\Pi_n$. Let $\vec{a}, b$ be sets. By strong $\Pi_n$-collection, there exists a set $C$ such that
$$(\forall x \in b)(\exists y \phi(y, x, \vec{a}) \Rightarrow (\exists y \in C)\phi(y, x, \vec{a})).$$
Therefore, using Lemma \ref{Th:Normalisation} and $\Pi_n$-separation,
$$A= \{x \in b \mid \exists y \phi(y, x, \vec{a})\}= \{x \in b \mid (\exists y \in C)\phi(y, x, \vec{a})\}$$
is a set. This completes the proof of (1).\\
We turn our attention to (2). Work in the theory $\mathbf{M}^-+ \Pi_n\textrm{-collection}+\Sigma_{n+1}\textrm{-separation}$. Let $\phi(x, y, \vec{z})$ be a $\Pi_n$-formula, and let $\vec{a}, b$ be sets. Now, $\Sigma_{n+1}$-separation implies that
$$A= \{x \in b \mid \exists y \phi(x, y, \vec{a})\}$$
is a set. And, $(\forall x \in A)\exists y \phi(x, y, \vec{a})$ holds. Therefore, we can apply $\Pi_n$-collection to obtain a set $C$ such that $(\forall x \in A)(\exists y \in C) \phi(x, y, \vec{a})$ holds. It now follows from the definition of $A$ that
$$(\forall x \in b)(\exists y \phi(x, y, \vec{a}) \Rightarrow (\exists y \in C) \phi(x, y, \vec{a})).$$
This completes the proof of (2).      
\Square
\end{proof}

\begin{Coroll1}
$\mathrm{MOST}$ ($\mathrm{MOST}^{-\mathrm{AC}}$, respectively) is the same theory as $\mathrm{Mac}+\textrm{strong } \Delta_0\textrm{-collection}$ ($\mathbf{M}+\textrm{strong } \Delta_0\textrm{-collection}$, respectively). \Square
\end{Coroll1}

Sufficiently rich set theories such as $\mathbf{M}$ and $\mathrm{KPI}$ allow us to express satisfaction in set structures. The following can be found in \cite{mat69} and \cite[Section III.1]{bar75}:

\begin{Lemma1} \label{Th:SatisfactionForSetStructuresInKPI}
In the theory $\mathrm{KPI}$, if $\mathcal{M}$ is a set structure, $\vec{a}$ is sequence of sets, and $\phi$ is an $\mathcal{L}(\mathcal{M})$-formula in the sense of the model whose arity agrees with the length of $\vec{a}$, then the predicate ``$\mathcal{M} \models \phi[\vec{v}/ \vec{a}]$" is definable by a $\Delta_1$-formula. \Square  
\end{Lemma1}

It is noted in \cite{mat01} that when powerset is present the recursions involved in the definition of satisfaction can be contained in sets even without any collection. The following is a consequence \cite[Proposition 3.10]{mat01}:

\begin{Lemma1}
In the theory $\mathbf{M}$, if $\mathcal{M}$ is a set structure, $\vec{a}$ is sequence of sets, and $\phi$ is an $\mathcal{L}(\mathcal{M})$-formula in the sense of the model whose arity agrees with the length of $\vec{a}$, then the predicate ``$\mathcal{M} \models \phi[\vec{v}/ \vec{a}]$" is definable and
$$\{\langle \ulcorner \phi \urcorner, \vec{a}\rangle \mid \vec{a} \in M \land \mathcal{M} \models \phi(\vec{a})\}$$
is a set.
\end{Lemma1}  

Equipped with these results, we can now define formulae that, in the theories $\mathrm{KPI}$ and $\mathbf{M}$, express satisfaction in the universe for the L\'{e}vy classes of $\mathcal{L}$-formulae.

\begin{Definitions1} \label{Df:Delta0Satisfaction}
Define $\mathrm{Sat}_{\Delta_0}(n, x)$ to be the formula
$$\begin{array}{c}
(n \in \omega) \land (n= \ulcorner \phi(v_1, \ldots, v_m) \urcorner \textrm{ where } \phi \textrm{ is } \Delta_0) \land (x= \langle x_1, \ldots, x_m \rangle) \land\\
\exists N \left( \bigcup N \subseteq N \land (x_1, \ldots, x_m \in N) \land (\langle N, \in \rangle \models \phi[x_1, \ldots, x_m]) \right) 
\end{array}.$$
\end{Definitions1}

The absoluteness of $\Delta_0$ properties between transitive structures and the universe, and the availability of $\mathrm{TCo}$ in $\mathrm{KPI}$ implies that the formula $\mathrm{Sat}_{\Delta_0}$ is equivalent, in the theory $\mathrm{KPI}$, to the formula
$$\begin{array}{c}
(n \in \omega) \land (n= \ulcorner \phi(v_1, \ldots, v_m) \urcorner \textrm{ where } \phi \textrm{ is } \Delta_0) \land (x= \langle x_1, \ldots, x_m \rangle) \land\\
\forall N \left( \bigcup N \subseteq N \land (x_1, \ldots, x_m \in N) \Rightarrow (\langle N, \in \rangle \models \phi[x_1, \ldots, x_m]) \right) 
\end{array}.$$
Therefore, Lemma \ref{Th:SatisfactionForSetStructuresInKPI} implies that $\mathrm{Sat}_{\Delta_0}(n, x)$ is $\Delta_1^{\mathrm{KPI}}$, and $\mathrm{Sat}_{\Delta_0}(n, x)$ expresses satisfaction for $\Delta_0$-formulae in the theories $\mathrm{KPI}$ and $\mathbf{M}$. We can now inductively define formulae $\mathrm{Sat}_{\Sigma_m}(n, x)$ and $\mathrm{Sat}_{\Pi_m}(n, x)$ that express satisfaction for formulae in the classes $\Sigma_m$ and $\Pi_m$.

\begin{Definitions1}
The formulae $\mathrm{Sat}_{\Sigma_m}(n, x)$ and $\mathrm{Sat}_{\Pi_m}(n, x)$ are defined inductively. Define $\mathrm{Sat}_{\Sigma_{m+1}}(n, x)$ to be the formula
$$\exists \vec{y} \exists k \exists b \left( \begin{array}{c}
(n= \ulcorner\exists \vec{u} \phi(\vec{u}, v_1, \ldots, v_l)\urcorner \textrm{ where } \phi \textrm{ is } \Pi_m)\land (x= \langle x_1, \ldots, x_l \rangle)\\
\land (b= \langle \vec{y}, x_1, \ldots, x_l \rangle) \land (k= \ulcorner \phi(\vec{u}, v_1, \ldots, v_l) \urcorner) \land \mathrm{Sat}_{\Pi_m}(k, b) 
\end{array}\right).$$
Define $\mathrm{Sat}_{\Pi_{m+1}}(n, x)$ to be the formula
$$\forall \vec{y} \forall k \forall b \left( \begin{array}{c}
(n= \ulcorner\forall \vec{u} \phi(\vec{u}, v_1, \ldots, v_l) \urcorner \textrm{ where } \phi \textrm{ is } \Sigma_m)\land (x= \langle x_1, \ldots, x_l \rangle)\\
\land ((b= \langle \vec{y}, x_1, \ldots, x_l \rangle) \land (k= \ulcorner\phi(\vec{u}, v_1, \ldots, v_l)\urcorner) \Rightarrow \mathrm{Sat}_{\Sigma_m}(k, b))
\end{array}\right).$$
\end{Definitions1}

The formula $\mathrm{Sat}_{\Sigma_m}(n, x)$ (respectively $\mathrm{Sat}_{\Pi_m}(n, x)$) is $\Sigma_m^{\mathrm{KPI}}$ ($\Pi_m^{\mathrm{KPI}}$, respectively), and, in the theories $\mathrm{KPI}$ and $\mathbf{M}$, expresses satisfaction for $\Sigma_m$-formulae ($\Pi_m$-formulae, respectively).

Another important feature of the theory $\mathrm{KPI}$ is its ability to construct $L$. The following can be found in \cite{mat69} and \cite[Chapter II]{bar75}:

\begin{Theorems1} \label{Th:DefinabilityOfLInKPI}
($\mathrm{KPI}$) The function $\alpha \mapsto L_\alpha$, where $\alpha$ is an ordinal, is total and $\Delta_1$. \Square 
\end{Theorems1}

As is usual, we use $V=L$ to abbreviate the expression that says that every set is the member of some $L_\alpha$ ($\forall x \exists \alpha((\alpha \textrm{ is an ordinal} \land (x \in L_\alpha)))$).     

We now turn to noting some of the properties of the theories $\mathbf{M}$, $\mathrm{Mac}$, $\mathbf{M}+H$ and $\mathrm{MOST}$ that are established in \cite{mat01}. The following useful fact is a consequence of \cite[Theorem I.6.1.]{bar75}:

\begin{Lemma1}
The theory $\mathrm{KPI}$ proves $\mathrm{TCo}$. \Square
\end{Lemma1} 

We also record the following consequence of \cite[Theorem Scheme 6.9(i)]{mat01}:

\begin{Theorems1} \label{Th:Delta0PSeparationInM}
The theory $\mathbf{M}$ proves all instances of $\Delta_0^\mathcal{P}$-separation. \Square
\end{Theorems1}

Section 2 of \cite{mat01} shows that by considering classes of well-founded extensional relations in a model of $\mathbf{M}$ one can obtain a model of $\mathbf{M}+\mathrm{H}$.

\begin{Theorems1} \label{Th:MisConsistentWithH}
(Mathias) If $\mathbf{M}$ is consistent, then so is $\mathbf{M}+\mathrm{H}$. 
\end{Theorems1}  

Section 3 of \cite{mat01} establishes a variety of consequences of Axiom H over the theories $\mathbf{M}$ and $\mathrm{Mac}$. A key observation of this section is that the theory $\mathrm{MOST}$ is exactly $\mathrm{Mac}$ plus Axiom H.

\begin{Lemma1} \label{Th:MOSTIsMacAndH}
$\mathrm{MOST}$ is the same theory as $\mathrm{Mac}+\textrm{Axiom }\mathrm{H}$. \Square
\end{Lemma1}

The following useful consequences of the theory $\mathrm{MOST}$ (=$\mathrm{Mac}+\textrm{strong } \Delta_0\textrm{-collection}$) are also proved in section 3 of \cite{mat01}:

\begin{Lemma1} \label{Th:BasicConsequencesMOST}
The theory $\mathrm{MOST}$ proves
\begin{itemize}
\item[(i)] every well-ordering is isomorphic to an ordinal,
\item[(ii)] every well-founded extensional relation is isomorphic to a transitive set,
\item[(iii)] for all cardinals $\kappa$, $\kappa^+$ exists,
\item[(iv)] for all cardinals $\kappa$, $H_\kappa$ exists.
\end{itemize}
\Square
\end{Lemma1}

Section 4 of \cite{mat01} establishes that the theory $\mathbf{M}+\mathrm{H}$ is capable of building G\"{o}del's $L$. Combined with Theorems \ref{Th:MisConsistentWithH} and \ref{Th:MOSTIsMacAndH} this yields that following consistency result:

\begin{Theorems1} \label{Th:MOSTplusVequalsLConsistentWithM}
(Mathias \cite[Theorem 1]{mat01}) If $\mathbf{M}$ is consistent, then so is $\mathrm{MOST}+V=L$. \Square 
\end{Theorems1}

The classes $\Delta_0^\mathcal{P}$, $\Sigma_1^\mathcal{P}$, $\Pi_1^\mathcal{P}$, \ldots are introduced and studied by Takahashi in \cite{tak72} where it is shown that for all $n \geq 1$, $(\Sigma_n^\mathcal{P})^{\mathrm{ZFC}}= \Sigma_{n+1}^\mathrm{ZFC}$, $(\Pi_n^\mathcal{P})^{\mathrm{ZFC}}= \Pi_{n+1}^\mathrm{ZFC}$, and $(\Delta_n^\mathcal{P})^{\mathrm{ZFC}}= \Delta_{n+1}^\mathrm{ZFC}$. The following calibration of Takahashi's result appears as Proposition Scheme 6.12 of \cite{mat01}:

\begin{Lemma1} \label{Th:TakahashiBaseCase}
(Takahashi) $\Sigma_1 \subseteq (\Delta_1^\mathcal{P})^{\mathrm{MOST}}$ and $\Delta_0^\mathcal{P} \subseteq \Delta_2^{\mathbf{S}_1}$. \Square
\end{Lemma1}

\noindent This yields the following refined version of Theorem 6 of \cite{tak72}:

\begin{Theorems1} \label{Th:TakahashiVsLevy}
(Takahashi) For all $n \geq 1$, $\Sigma_{n+1} \subseteq (\Sigma_n^\mathcal{P})^{\mathrm{MOST}}$, $\Pi_{n+1} \subseteq (\Pi_n^\mathcal{P})^{\mathrm{MOST}}$, $\Delta_{n+1} \subseteq (\Delta_n^\mathcal{P})^{\mathrm{MOST}}$, $\Sigma_n^\mathcal{P} \subseteq \Sigma_{n+1}^{\mathbf{S}_1}$, $\Pi_n^\mathcal{P} \subseteq \Pi_{n+1}^{\mathbf{S}_1}$, and $\Delta_n^\mathcal{P} \subseteq \Delta_{n+1}^{\mathbf{S}_1}$. \Square
\end{Theorems1} 

Lemmas \ref{Th:BasicRelationships} and \ref{Th:StrongCollectionVsCollectionSeparation}, and Theorem \ref{Th:TakahashiVsLevy} now show:

\begin{Coroll1}
The theory $\mathbf{M}+\textrm{strong } \Pi_1\textrm{-collection}$ proves every axiom of $\mathrm{KP}^\mathcal{P}$. \Square
\end{Coroll1}

In \cite{mat01}, Mathias proves a $\Sigma_1^\mathcal{P}$-Recursion Theorem in the theory $\mathrm{KP}^\mathcal{P}$. The following appear as Lemma 6.25 and Theorem 6.26 in \cite{mat01}:

\begin{Lemma1}
If $F$ is a total $\Sigma_1^\mathcal{P}$-definable class function, then the formula $y=F(x)$ is $\Delta_1^\mathcal{P}$. \Square 
\end{Lemma1}

\begin{Theorems1} \label{Th:Sigma1PRecursionTheorem}
($\mathrm{KP}^\mathcal{P}$) Let $G$ be a $\Sigma_1^\mathcal{P}$-definable class. If $G$ is a total function, then there exists a $\Sigma_1^\mathcal{P}$-definable total class function $F$ such that for all $x$, $F(x)= G(F\upharpoonright x)$. \Square
\end{Theorems1}

The fact that we have access to Theorem \ref{Th:Sigma1PRecursionTheorem} in the theory $\mathbf{M}+\textrm{strong } \Pi_1\textrm{-collection}$ yields:

\begin{Coroll1}
The theory $\mathbf{M}+\textrm{strong } \Pi_1\textrm{-collection}$ proves that for all ordinals $\alpha$, $V_\alpha$ is a set. Moreover, the formula ``$x= V_\alpha$" with free variables $x$ and $\alpha$ is equivalent to a $\Delta_1^\mathcal{P}$-formula. \Square
\end{Coroll1}

Results proved in \cite{mat01} also reveal that the theory $\mathbf{M}+\textrm{strong } \Pi_1\textrm{-collection}$ is capable of proving the consistency of Zermelo Set Theory plus $\Delta_0$-collection. Mathias \cite[Lemma 6.31]{mat01} shows that the theory obtained by strengthening $\mathrm{KP}$ with an axiom that asserts the existence of $V_\alpha$ for every ordinal $\alpha$ is capable of proving the consistency of $\mathrm{Z}$. The fact that $\mathrm{KP}^\mathcal{P}$ is equipped with enough recursion to prove the existence of $V_\alpha$ for every $\alpha$ \cite[Proposition 6.28]{mat01} thus yields:    

\begin{Theorems1} \label{Th:KPPProvesConsistencyOfZ}
(Mathias) The theory $\mathrm{KP}^\mathcal{P}$ proves that there exists a transitive model of $\mathrm{Z}$.\Square
\end{Theorems1}

Mathias \cite[Theorem 5]{mat01} also shows that all of the axioms of $\mathrm{KP}$ plus $V=L$ can be consistently added to $\mathrm{Z}$. In particular:

\begin{Theorems1} \label{Th:ZisEquiconsistentWithKLZ}
(Mathias) If $\mathrm{Z}$ is consistent, then so is $\mathrm{Z}+\Delta_0\textrm{-collection}+V=L$. \Square   
\end{Theorems1}

Theorems \ref{Th:KPPProvesConsistencyOfZ} and \ref{Th:ZisEquiconsistentWithKLZ} now yield:

\begin{Coroll1}
$\mathrm{KP}^\mathcal{P} \vdash \mathrm{Con}(\mathrm{Z}+\Delta_0\textrm{-collection}+V=L)$. \Square
\end{Coroll1}

\section[The strength of $\Delta_0^\mathcal{P}$-collection]{The strength of $\Delta_0^\mathcal{P}$-collection}\label{Sec:BaseCaseSection}

In this section we investigate the strength of adding $\Delta_0^\mathcal{P}$-collection to subsystems of set theory studied in \cite{mat01}. We show that if $T$ is one of the theories $\mathbf{M}$, $\mathbf{M}+\mathrm{H}$, $\mathrm{Mac}$ or $\mathrm{MOST}$, then the theory obtained by adding $\Delta_0^\mathcal{P}$-collection to $T$ is $\Pi_2^\mathcal{P}$-conservative over $T$. Combined with Theorems \ref{Th:MOSTplusVequalsLConsistentWithM} and \ref{Th:TakahashiVsLevy}, this shows that if $\mathbf{M}$ is consistent, then so is $\mathrm{MOST}+\Pi_1\textrm{-collection}$. 

If $u$ is a set, then we will $H_{\leq|u|}$ to denote the set 
$$\{x \mid |\mathrm{TC}(\{x\})|\leq |u|\}.$$

\begin{Lemma1}\label{Th:HuExist}
The theory $\mathbf{M}+\mathrm{H}$ proves that for all sets $u$, $H_{\leq|u|}$ exists.
\end{Lemma1}

\begin{proof}
Work in the theory $\mathbf{M}+\mathrm{H}$. Let $u$ be a set. Using Axiom H, let $T$ be a set such that
$$\forall z\left(\bigcup z \subseteq z \land |z|\leq |u| \Rightarrow z \subseteq T\right)$$
Note that if $x$ is a set such that $|\mathrm{TC}(\{x\})| \leq |u|$, then $\mathrm{TC}(\{x\}) \subseteq T$ and so $x \in T$. Moreover, if $|\mathrm{TC}(\{x\})| \leq |u|$, then $\mathrm{TC}(\{x\}) \in \mathcal{P}(T)$ and the injection witnessing $|\mathrm{TC}(\{x\})| \leq |u|$ is in $\mathcal{P}(T\times u)$. Therefore $\Delta_0$-separation implies that $H_{\leq|u|}$ exists.
\Square   
\end{proof}

The following is immediate from the definition of $H_{\leq|u|}$:

\begin{Lemma1}\label{Th:HuSuperTransitive} 
The theory $\mathbf{M}+\mathrm{H}$ proves that if $u, x, y$ are sets, then
\begin{itemize}
\item[(i)] if $x \in y \in H_{\leq|u|}$, then $x \in H_{\leq|u|}$, and
\item[(ii)] if $x \subseteq y \in H_{\leq|u|}$, then $x \in H_{\leq|u|}$. 
\end{itemize}
\Square
\end{Lemma1}

\begin{Definitions1} \label{Df:HApproximations}
Let $n \in \omega$ and let $u$ be a set. We say that $f$ is an $n$-good $|u|$-$H$-approximation if 
\begin{itemize}
\item[(i)] $f$ is a function and $\mathrm{dom}(f)= n+1$
\item[(ii)] $f(\emptyset)= H_{\leq|u|}$ 
\item[(iii)] $(\forall k \in n+1)\exists v (f(k)= H_{\leq|v|})$
\item[(iv)] $(\forall k \in n)(f(k) \in f(k+1))$.
\end{itemize}
\end{Definitions1}

We first observe that in any model of $\mathbf{M}+\mathrm{H}$ there exists an $n$-good $|u|$-$H$-approximation for every externally finite $n$ and every set $u$ in the model.

\begin{Lemma1} \label{Th:FiniteHApproximationsExist}
Let $n \in \omega$. If $\mathcal{M} \models \mathbf{M}+\mathrm{H}$ and $u \in M$, then 
$$\mathcal{M} \models \exists f(f \textrm{ is an } n\textrm{-good } |u|\textrm{-}H\textrm{-approximation}).$$
\end{Lemma1}

\begin{proof}
Let $\mathcal{M}= \langle M, \in^\mathcal{M} \rangle$ be such that $\mathcal{M} \models \mathbf{M}+\mathrm{H}$ and let $u \in M$. We prove, by external induction on $\omega$, that for all $n \in \omega$,
$$\mathcal{M} \models \exists f(f \textrm{ is an } n\textrm{-good } |u|\textrm{-}H\textrm{-approximation}).$$
It follows from Lemma \ref{Th:HuExist} that
$$\mathcal{M} \models \exists f(f \textrm{ is a }0 \textrm{-good } |u|\textrm{-}H\textrm{-approximation}).$$
Suppose that the lemma is false, and $k \in \omega$ is least such that
$$\mathcal{M} \models \neg \exists f(f \textrm{ is a } (k+1)\textrm{-good } |u|\textrm{-}H\textrm{-approximation}).$$
Work inside $\mathcal{M}$. Let $f$ be a $k$-good $|u|$-$H$-approximation. Let $v= f(k)\cup\{f(k)\}$. It follows from Definition \ref{Df:HApproximations}(iii) and Lemma \ref{Th:HuSuperTransitive} that $v= \mathrm{TC}(\{f(k)\})$. Therefore $g=f\cup \{\langle k+1, H_{\leq|v|}\rangle\}$ is a $(k+1)$-good $|u|$-$H$-approximation, which is a contradiction.   
\Square
\end{proof}

In the proof of the following result we obtain models of $\Delta_0^\mathcal{P}$-collection by considering a cut of an $n$-good $|u|$-$H$-approximation of nonstandard length. This idea of obtaining ``more" collection from a cut of a nonstandard model of set theory also appears in Ressayre's work on limitations of extensions of Kripke-Platek Set Theory \cite{res87} (see also \cite{flw16}) and Friedman's work \cite{fri73} on the standard part of countable non-standard models of set theory.

\begin{Theorems1}\label{Th:Delta0PCollectionConsistentWithMOST}
\begin{itemize}
\item[(I)] The theory $\mathbf{M}+\mathrm{H}+\Delta_0^\mathcal{P}\textrm{-collection}$ is $\Pi_2^\mathcal{P}$-conservative over the theory $\mathbf{M}+\mathrm{H}$.
\item[(II)] The theory $\mathrm{MOST}+\Pi_1\textrm{-collection}$ is $\Pi_3$-conservative over the theory $\mathrm{MOST}$. 
\end{itemize}
\end{Theorems1}

\begin{proof}
To prove (I) it is sufficient to show that every $\Sigma_2^\mathcal{P}$-sentence that is consistent with $\mathbf{M}+\mathrm{H}$ is also consistent with $\mathbf{M}+\mathrm{H}+\Delta_0^\mathcal{P}\textrm{-collection}$. Suppose that $\exists \vec{x} \forall \vec{y} \theta(\vec{x}, \vec{y})$, where $\theta(\vec{x}, \vec{y})$ is a $\Delta_0^\mathcal{P}$-formulae, is consistent with $\mathbf{M}+\mathrm{H}$. Let $\mathcal{M}= \langle M, \in^\mathcal{M} \rangle$ be a recursively saturated model of $\mathbf{M}+\mathrm{H}+\exists \vec{x} \forall \vec{y} \theta(\vec{x}, \vec{y})$. Let $\vec{a} \in M$ be such $\mathcal{M}\models \forall \vec{y} \theta(\vec{a}, \vec{y})$ and let $u \in M$ be such that $\vec{a} \in u$. Consider the type
$$\Xi(x, u)= \{x \in \omega\}\cup\{x > n \mid n \in \omega\}\cup \{\exists f(f \textrm{ is an }x\textrm{-good } u\textrm{-}H\textrm{-approximation})\}.$$
By Lemma \ref{Th:FiniteHApproximationsExist}, $\Xi(x, u)$ is finitely realised in any model of $\mathbf{M}+\textrm{H}$, and so there exists $k \in M$ such that $\Xi(k, u)$ is satisfied in $\mathcal{M}$. Note that $k$ is a nonstandard element of $\omega^{\mathcal{M}}$. Let $f \in M$ be such that 
$$\mathcal{M} \models (f \textrm{ is a }k\textrm{-good } u\textrm{-}H\textrm{-approximation}).$$
Define $\mathcal{N}= \langle N, \in^\mathcal{N} \rangle$ by
$$N= \bigcup_{n \in \omega} f(n^\mathcal{M})^* \textrm{ and } \in^\mathcal{N} \textrm{ is the restriction of } \in^\mathcal{M} \textrm{ to }N.$$
We claim that $\mathcal{N}$ satisfies $\mathbf{M}+\textrm{H}+\Delta_0^\mathcal{P}\textrm{-collection}+\exists \vec{x} \forall \vec{y} \theta(\vec{x}, \vec{y})$. Note that $\mathcal{N} \subseteq_e^\mathcal{P} \mathcal{M}$ and $\vec{a} \in N$, so $\mathcal{N}\models \exists \vec{x} \forall \vec{y} \theta(\vec{x}, \vec{y})$. Let $x \in N$. Let $n \in \omega$ be such that $\mathcal{M} \models (x \in f(n^\mathcal{M}))$. Therefore $\mathcal{M}\models (\mathcal{P}(x) \subseteq f(n^\mathcal{M}))$ and $f(n^{\mathcal{M}}) \in (f((n+1)^\mathcal{M}))^* \subseteq N$. It now follows from Definition \ref{Df:HApproximations} that $\mathcal{P}^\mathcal{M}(x) \in N$. Therefore $\mathcal{N}\models (\mathrm{powerset})$ and for all $x \in N$, $\mathcal{P}^\mathcal{N}(x)= \mathcal{P}^\mathcal{M}(x)$. It is now clear that $\mathcal{N} \models \mathbf{M}$. 

We turn to showing that Axiom H holds in $\mathcal{N}$. Let $u \in N$. Let $n \in \omega$ be such that $u \in f(n^\mathcal{M})^*$. By Definition \ref{Df:HApproximations}, there exists $v \in M$ such that $\mathcal{M} \models (f(n^\mathcal{M})= H_{\leq|v|})$, and so $\mathcal{M} \models (|u|\leq |v|)$. Now, working inside $\mathcal{N}$, if $z$ is transitive with $|z| \leq |u|$, then $|z|\leq |v|$ and so $z \in f(n^{\mathcal{M}})$. Therefore
$$\mathcal{N} \models \forall z\left(\bigcup z \subseteq z \land |z|\leq |u| \Rightarrow z \in f(n^{\mathcal{M}})\right)$$
and so Axiom H holds in $\mathcal{N}$. 

We are left to show that $\mathcal{N}$ satisfies $\Delta_0^\mathcal{P}\textrm{-collection}$. We make use of the following property of $\mathcal{N}$:\\
{\bf Claim:} If $C\in M$ and $C^* \subseteq N$, then $C \in N$.\\
We prove this claim. Suppose, for a contradiction, that $C \in M$, $C^* \subseteq N$ and $C \notin N$. Note that if $n \in k^*$ is nonstandard, then $C^* \subseteq f(n)^*$ and $\mathcal{M} \models (C \in f(n+1))$. Therefore, working inside $\mathcal{M}$, the set 
$$A= \{n \in k \mid C \notin f(n)\}$$
defines the standard $\omega$, which is a contradiction. This proves the claim.\\
Now, let $\phi(x, y, \vec{z})$ be a $\Delta_0^\mathcal{P}$-formula. Let $\vec{d}, b \in N$ be such that
$$\mathcal{N}\models (\forall x \in b) \exists y \phi(x, y, \vec{d}).$$
The following formula is a $\Delta_0^\mathcal{P}$-formula with parameters $\vec{d}$, $k$ and $f$:
$$\phi(x, y, \vec{d}) \land (\forall n \in k)(y \notin f(n) \Rightarrow \neg (\exists w \in f(n))\phi(x, w, \vec{d})).$$
So, by $\Delta_0^\mathcal{P}$-absoluteness,
$$\mathcal{M} \models (\forall x \in b)(\exists y \in f(k))(\phi(x, y, \vec{d}) \land (\forall n \in k)(y \notin f(n) \Rightarrow \neg (\exists w \in f(n))\phi(x, w, \vec{d}))).$$
Working inside $\mathcal{M}$, $\Delta_0^\mathcal{P}$-separation (Theorem \ref{Th:Delta0PSeparationInM}) implies that
$$C= \{\langle x, y \rangle \in b \times f(k) \mid \phi(x, y, \vec{d}) \land (\forall n \in k)(y \notin f(n) \Rightarrow \neg (\exists w \in f(n))\phi(x, w, \vec{d}))$$
is a set. And $\Delta_0^\mathcal{P}$-absoluteness implies that $C^* \subseteq N$. Therefore $C \in N$. Working inside $\mathcal{N}$, let $B= \mathrm{rng}(C)$. So, 
$$\mathcal{N} \models (\forall x \in b)(\exists y \in B) \phi(x, y, \vec{d}),$$
which shows that $\mathcal{N} \models \Delta_0^\mathcal{P}\textrm{-collection}$.

To see that (II) holds observe that if the Axiom of Choice holds in $\mathcal{M}$ in the proof of (I), then it also holds in $\mathcal{N}$. It then follows from Theorem \ref{Th:TakahashiVsLevy} that $\mathcal{N}$ also satisfies $\Pi_1$-collection, and we get $\Pi_3$-conservativity.              
\Square
\end{proof}

Theorem \ref{Th:Delta0PCollectionConsistentWithMOST} combined with Theorems \ref{Th:MOSTplusVequalsLConsistentWithM} shows that the consistency $\mathbf{M}$ implies the consistency of $\mathrm{MOST}+\Pi_1\textrm{-collection}$.

\begin{Coroll1}\label{Th:ExtensionOfEquiconsistencyOfM}
If $\mathbf{M}$ is consistent, then so is $\mathrm{MOST}+\Pi_1\textrm{-collection}$ ($= \mathrm{Mac}+\Pi_1\textrm{-collection}$). \Square 
\end{Coroll1}

The argument used in the proof of Theorem \ref{Th:Delta0PCollectionConsistentWithMOST} can also be used to show that that the theories $\mathbf{M}+\Delta_0^\mathcal{P}\textrm{-collection}$ and $\mathrm{Mac}+\Delta_0^\mathcal{P}\textrm{-collection}$ are $\Pi_2^\mathcal{P}$-conservative over the theories $\mathbf{M}$ and $\mathrm{Mac}$, respectively. To see this we introduce a modification of Definition \ref{Df:HApproximations}:

\begin{Definitions1}\label{Df:PowersetApproximations}
Let $n \in \omega$ and let $u$ be a set. We say that $f$ is an $n$-good $u$-$\mathcal{P}$-approximation if 
\begin{itemize}
\item[(i)] $f$ is a function and $\mathrm{dom}(f)= n+1$
\item[(ii)] $f(\emptyset)= \mathrm{TC}(u)$ 
\item[(iii)] $(\forall k \in n)(f(k+1)= \mathcal{P}(f(k)))$.
\end{itemize}
\end{Definitions1}  

An $n$-good $u$-$\mathcal{P}$-approximation is a sequence $\mathcal{P}(v)$, $\mathcal{P}(\mathcal{P}(v)), \ldots$ where $v$ is the transitive closure of $u$. The same argument that was used to prove Lemma \ref{Th:FiniteHApproximationsExist} shows that in any model of $\mathbf{M}$, any such sequence with externally finite length is guaranteed to exist.

\begin{Lemma1}\label{Th:FinitePowersetApproximationsExist}
Let $n \in \omega$. If $\mathcal{M} \models \mathbf{M}$ and $u \in M$, then 
$$\mathcal{M} \models \exists f(f \textrm{ is an } n\textrm{-good } u\textrm{-}\mathcal{P}\textrm{-approximation}).$$
\Square
\end{Lemma1}

Replacing the $n$-good $|u|$-$H$-approximations in the proof of Theorem \ref{Th:Delta0PCollectionConsistentWithMOST} now shows that adding $\Delta_0^\mathcal{P}$-collection to $\mathbf{M}$ or $\mathrm{Mac}$ does not prove any new $\Pi_2^\mathcal{P}$-sentences.

\begin{Theorems1}\label{Th:Delta0PCollectionConsistentWithMac}
\begin{itemize}
\item[(I)] The theory $\mathbf{M}+\Delta_0^\mathcal{P}\textrm{-collection}$ is $\Pi_2^\mathcal{P}$-conservative over the theory $\mathbf{M}$.
\item[(II)] The theory $\mathrm{Mac}+\Delta_0^\mathcal{P}\textrm{-collection}$ is $\Pi_2^\mathcal{P}$-conservative over the theory $\mathrm{Mac}$. 
\end{itemize}
\Square
\end{Theorems1} 

\begin{Remark1}
Theorems \ref{Th:Delta0PCollectionConsistentWithMOST} and \ref{Th:Delta0PCollectionConsistentWithMac} highlight a mistake in the final sentence of \cite[Metatheorem 9.41]{mat01} and the final clause, starting after the colon, of \cite[Theorem 16]{mat01} (which paraphrases \cite[Metatheorem 9.41]{mat01}). This erroneous assertion is used by the author in \cite{mck15} to claim that the theory $\mathrm{Mac}+\Delta_0^\mathcal{P}\textrm{-collection}$ represents a new lower-bound on the consistency strength of the theory $\mathrm{NFU}+\mathrm{AxCount}_\leq$. Theorem \ref{Th:Delta0PCollectionConsistentWithMOST} now shows that $\mathrm{Mac}+\Delta_0^\mathcal{P}\textrm{-collection}$ does not represent an improvement on previously known lower-bounds on the consistency strength of $\mathrm{NFU}+\mathrm{AxCount}_\leq$. 
\end{Remark1}

\section[The strength of $\Pi_n$-collection over $\mathbf{M}$]{The strength of $\Pi_n$-collection over $\mathbf{M}$} \label{Sec:GeneralCollectionResults}

In this section we generalise and expand upon Theorem \ref{Th:Delta0PCollectionConsistentWithMOST} to show for all $n \geq 1$,
\begin{enumerate}
\item the theory $\mathbf{M}+\Pi_{n+1}\textrm{-collection}$ is $\Pi_{n+3}$-conservative over the theory\\ 
$\mathbf{M}+\textrm{strong }\Pi_n\textrm{-collection}$,
\item the theory $\mathbf{M}+\Pi_{n+1}\textrm{-collection}+ \Sigma_{n+2}\textrm{-induction on } \omega$ proves that there exists a transitive model of $\mathrm{Z}+ \Pi_n\textrm{-collection}$.
\end{enumerate}
The main tool used in the proof of these results will be the following modification and generalisation of Definition \ref{Df:HApproximations}:

\begin{Definitions1} \label{Df:nGoodSubmodelApproximation}
Let $n, m \in \omega$, and let $\alpha$ be an ordinal. We say that $f$ is an $n$-good $\langle m+1, \alpha \rangle$-submodel approximation if 
\begin{itemize}
\item[(i)] $f$ is a function and $\mathrm{dom}(f)= n+1$ 
\item[(ii)] $f(\emptyset)= V_\alpha$
\item[(iii)] $(\forall k \in n+1)\exists \beta((\beta \textrm{ is an ordinal})\land f(k)= V_\beta)$
\item[(iv)]
$$(\forall k \in n)(\forall l \in \omega)(\forall a \in f(k+1))((\langle f(k+1), \in\rangle \models \mathrm{Sat}_{\Pi_m}(l, a)) \Rightarrow \mathrm{Sat}_{\Pi_m}(l, a))$$
\item[(v)]
$$(\forall k \in n)(\forall l \in \omega)(\forall a \in f(k))(\mathrm{Sat}_{\Sigma_{m+1}}(l, a) \Rightarrow (\langle f(k+1), \in  \rangle \models \mathrm{Sat}_{\Sigma_{m+1}}(l, a)))$$
\end{itemize}
\end{Definitions1}

An $n$-good $\langle m+1, \alpha\rangle$-submodel approximation is a sequence $\langle V_{\beta_0}, \ldots, V_{\beta_n} \rangle$ such that $V_{\beta_0}= V_\alpha$ (condition (ii)), for all $0 \leq l < k$, $\beta_l \leq \beta_k$ (condition (v) applied to the $\Sigma_1$-formula ``$\exists v(a \in v)$"), each $V_{\beta_k}$ ($1 \leq k \leq n$) is a $\Pi_m$-elementary submodel of the universe (condition (iv)), each $V_{\beta_{k+1}}$ satisfies the same $\Sigma_{m+1}$-formulae with parameters from $V_{\beta_k}$ as the universe (condition (v)). Note that if an infinite sequence $\langle V_{\beta_0}, V_{\beta_1}, \ldots \rangle$ is such that for every $n \in \omega$, the first $n+1$ elements of this sequence form an $n$-good $\langle m+1, \alpha \rangle$-submodel approximation, then $\bigcup_{n \in \omega} V_{\beta_n}$ is a $\Pi_{m+1}$-elementary submodel of the universe. 

We make the following observations about the complexity of Definition \ref{Df:nGoodSubmodelApproximation}:
\begin{enumerate}
\item The formula ``$f$ is a function and $\mathrm{dom}(f)= n+1$" is $\Delta_0$ with parameters $f$ and $n$.
\item The formula ``$f(\emptyset)= V_\alpha$" is $\Delta_0$ with parameters $f$ and $V_\alpha$. 
\item The formula ``$(\forall k \in n+1)\exists \beta((\beta \textrm{ is an ordinal})\land f(k)= V_\beta)$" is both $\Sigma_2^{\mathbf{M}+\textrm{strong } \Pi_1\textrm{-collection}}$ and $(\Sigma_1^\mathcal{P})^{\mathbf{M}+\textrm{strong } \Pi_1\textrm{-collection}}$ with parameters $f$ and $n$.
\item For all $m \in \omega$, the formula
$$(\forall k \in n)(\forall l \in \omega)(\forall a \in f(k+1))((\langle f(k+1), \in\rangle \models \mathrm{Sat}_{\Pi_m}(l, a)) \Rightarrow \mathrm{Sat}_{\Pi_m}(l, a))$$
is $\Pi_{\max(1, m)}^{\mathrm{KPI}}$ with parameters $f$ and $n$. 
\item For all $m \in \omega$, the formula
$$(\forall k \in n)(\forall l \in \omega)(\forall a \in f(k))(\mathrm{Sat}_{\Sigma_{m+1}}(l, a) \Rightarrow (\langle f(k+1), \in  \rangle \models \mathrm{Sat}_{\Sigma_{m+1}}(l, a)))$$
is $\Pi_{m+1}^{\mathrm{KPI}}$ with parameters $f$ and $n$.  
\end{enumerate}
In light of these observations we introduce specific notion for the formulae that say that $f$ is an $n$-good $\langle m+1, \alpha \rangle$-submodel approximation.

\begin{Definitions1}
Let $\alpha$ be an ordinal and let $m \in \omega$. We write $\Psi_m(n, f, V_\alpha)$ for the formula, with free variables $f$ and $n$, and parameter $V_\alpha$, that the theory $\mathbf{M}+\textrm{strong } \Pi_1\textrm{-collection}$ proves asserts that $f$ in an $n$-good $\langle m+1, \alpha\rangle$-submodel approximation, and such that $\Psi_0(n, f, V_\alpha)$ is $\Sigma_2$, $\Psi_1(n, f, V_\alpha)$ is $\mathbf{Bol}(\Sigma_2)$, and if $m > 1$, $\Psi_m(n, f, V_\alpha)$ is $\Pi_{m+1}$. 
\end{Definitions1}

\begin{Lemma1} \label{Th:BaseCaseIKeyLemma}
The theory $\mathbf{M}+\textrm{strong } \Pi_{1}\textrm{-collection}$ proves that for all ordinals $\alpha$ and for all $n \in \omega$, there exists an $n$-good $\langle 1, \alpha\rangle$-submodel approximation.
\end{Lemma1}

\begin{proof}
Work in the theory $\mathbf{M}+\textrm{strong } \Pi_{1}\textrm{-collection}$. Let $\alpha$ be an ordinal. We will use $\Sigma_2$-induction on $\omega$ to prove $(\forall n \in \omega) \exists f \Psi_0(n, f, V_\alpha)$. It is clear that $\exists f \Psi_0(\emptyset, f, V_\alpha)$ holds. Let $n \in \omega$ and suppose that $f$ is such that $\Psi_0(n, f, V_\alpha)$ holds. Let $\beta$ be the ordinal such that $f(n)= V_\beta$. Consider the $\Sigma_1$-formula $\psi(x, y)$ defined by
$$\exists z \exists a \exists l((x= \langle a, l \rangle) \land (z= \langle y, a \rangle) \land (l= \ulcorner \phi(u, v) \urcorner \textrm{ where } \phi \textrm{ is } \Delta_0)\land \mathrm{Sat}_{\Delta_0}(l, z)).$$
Strong $\Sigma_1$-collection implies that there exists a $C$ such that
$$(\forall x \in V_\beta \times \omega)(\exists y \psi(x, y) \Rightarrow (\exists y \in C) \psi(x, y)).$$
Let $\gamma > \beta$ be such that $C \subseteq V_\gamma$. Therefore, for all $l \in \omega$ and for all $a \in V_\beta$,
$$\textrm{if } \mathrm{Sat}_{\Sigma_1}(l, a), \textrm{ then } \langle V_\gamma, \in \rangle \models \mathrm{Sat}_{\Sigma_1}(l, a).$$
It now follows that $g= f \cup \{\langle n+1, V_\gamma \rangle\}$ satisfies $\Psi_0(n+1, g, V_\alpha)$. The fact that $(\forall n \in \omega) \exists f \Psi_0(n, f, V_\alpha)$ holds now follows by $\Sigma_2$-induction on $\omega$. 
\Square
\end{proof}

\begin{Lemma1} \label{Th:BaseCaseIIKeyLemma}
The theory $\mathbf{M}+\textrm{strong } \Pi_{1}\textrm{-collection}$ proves that for all ordinals $\alpha$, there exists a function $f$ with $\mathrm{dom}(f)=\omega$ such that for all $n \in \omega$, $f \upharpoonright (n+1)$ is an $n$-good $\langle 1, \alpha\rangle$-submodel approximation.
\end{Lemma1}

\begin{proof}
Work in the theory $\mathbf{M}+\textrm{strong } \Pi_{1}\textrm{-collection}$. Using Lemma \ref{Th:BaseCaseIKeyLemma} and strong $\Sigma_2$-collection, we can find a set $B$ such that $(\forall n \in \omega)(\exists f \in B) \Psi_0(n, f, V_\alpha)$ holds. Now, $\Sigma_2$-separation ensures that 
$$D= \{f \in B \mid (\exists n \in \omega)\Psi_0(n, f, V_\alpha)\}$$
is a set. Let
$$G= \left\{f \in D \Big| (\forall k \in \mathrm{dom}(f))(\forall g \in D)\left(\begin{array}{c}
(k \in \mathrm{dom}(g)) \land (g(k) \neq f(k))\\
\Rightarrow f(k) \in g(k)
\end{array} \right)\right\},$$
which is a set. Now, for all $f_1, f_2 \in G$, $f_1$ and $f_2$ agree on their common domain. Moreover, a straightforward internal induction using the fact that Lemma \ref{Th:BaseCaseIKeyLemma} holds shows that for all $n \in \omega$, $(\exists f \in G)(\mathrm{dom}(f)=n+1)$ holds. Therefore $g= \bigcup G$ is a function with domain $\omega$ such that for all $n \in \omega$, $\Psi_0(n, g \upharpoonright (n+1), V_\alpha)$ holds.
\Square
\end{proof}

We can now prove analogues of Lemmas \ref{Th:BaseCaseIKeyLemma} and \ref{Th:BaseCaseIIKeyLemma} for the theories $\mathbf{M}+ \Pi_{m}\textrm{-collection}+\Sigma_{m+1}\textrm{-induction on } \omega$ where $m \geq 2$.

\begin{Lemma1} \label{Th:KeyLemmaForProvingConsistency}
Let $m \geq 1$. The theory $\mathbf{M}+ \Pi_{m+1}\textrm{-collection}+\Sigma_{m+2}\textrm{-induction on } \omega$ proves
\begin{itemize}
\item[(I)] for all ordinals $\alpha$ and for all $n \in \omega$, there exists an $n$-good $\langle m+1, \alpha\rangle$-submodel approximation,
\item[(II)] for all ordinals $\alpha$, there exists a function $f$ with $\mathrm{dom}(f)=\omega$ such that for all $n \in \omega$, $f \upharpoonright (n+1)$ is an $n$-good $\langle m+1, \alpha\rangle$-submodel approximation.
\end{itemize}
\end{Lemma1}

\begin{proof}
We prove this lemma by external induction on $m$. We begin by proving the induction step. Suppose that (I) and (II) of the lemma hold for $m=p \geq 1$. Work in the theory $\mathbf{M}+\Pi_{p+2}\textrm{-collection}+\Sigma_{p+3}\textrm{-induction on } \omega$. Let $\alpha$ be an ordinal. We will use $\Sigma_{p+3}$-induction on $\omega$ to show that $(\forall n \in \omega) \exists f \Psi_{p+1}(n, f, V_\alpha)$ holds. It is clear that $\exists f \Psi_{p+1}(\emptyset, f, V_\alpha)$ holds. Let $n \in \omega$, and suppose that $\exists f \Psi_{p+1}(n, f, V_\alpha)$ holds. Let $f$ be such that $\Psi_{p+1}(n, f, V_\alpha)$. Let $\delta$ be the ordinal such that $f(n)= V_\delta$. Consider the $\Sigma_{p+2}$-formula $\psi(x, y)$ defined by
$$\exists z \exists a \exists l((x= \langle a, l \rangle) \land (z= \langle y, a \rangle) \land (l= \ulcorner \phi(u, v) \urcorner \textrm{ where } \phi \textrm{ is } \Pi_{p+1})\land \mathrm{Sat}_{\Pi_{p+1}}(l, z)).$$
Strong $\Sigma_{p+2}$-collection implies that there exists a $C$ such that
$$(\forall x \in V_\delta \times \omega)(\exists y \psi(x, y) \Rightarrow (\exists y \in C) \psi(x, y)).$$
Let $\beta > \delta$ be such that $C \subseteq V_\beta$. Now, using (II) of the induction hypothesis, we can find a function $g$ with $\mathrm{dom}(g)=\omega$ such that for all $q \in \omega$, $\Psi_p(q, g \upharpoonright (q+1), V_\beta)$. Now, let $\gamma > \beta$ be such that $V_\gamma= \bigcup \mathrm{rng}(g)$. It follows from (iv) and (v) of Definition (\ref{Df:nGoodSubmodelApproximation}) that for all $l \in \omega$ and for all $a \in V_\gamma$,
$$\textrm{if } \langle V_\gamma, \in \rangle \models \mathrm{Sat}_{\Pi_{p+1}}(l, a), \textrm{ then } \mathrm{Sat}_{\Pi_{p+1}}(l, a).$$
And, since $C \subseteq V_\beta \subseteq V_\gamma$, for all $l \in \omega$ and for all $a \in V_\delta$,
$$\textrm{if } \mathrm{Sat}_{\Sigma_{p+2}}(l, a), \textrm{ then } \langle V_\gamma, \in \rangle \models \mathrm{Sat}_{\Sigma_{p+2}}(l, a).$$
Therefore, the function $h=f \cup \{\langle n+1, V_\gamma \rangle\}$ satisfies $\Psi_{p+1}(n+1, h, V_\alpha)$. The fact that $(\forall n \in \omega)\exists f \Psi_{p+1}(n, f, V_\alpha)$ now follows from $\Sigma_{p+3}$-induction on $\omega$. This completes the induction step for (I). Turning our attention to (II), we can use $\Pi_{p+2}$-collection to find a set $B$ such that $(\forall n \in \omega)(\exists f \in B) \Psi_{p+1}(n, f, V_\alpha)$. Now, $\Pi_{p+2}$-separation ensures that
$$D= \{f \in B \mid (\exists n \in \omega) \Psi_{p+1}(n, f, V_\alpha)\}$$
is a set. Let
$$G= \left\{f \in D \Big| (\forall k \in \mathrm{dom}(f))(\forall g \in D)\left(\begin{array}{c}
(k \in \mathrm{dom}(g)) \land (g(k) \neq f(k))\\
\Rightarrow f(k) \in g(k)
\end{array} \right)\right\},$$
As in the proof of Lemma \ref{Th:BaseCaseIIKeyLemma}, if $f_1, f_2 \in G$, then $f_1$ and $f_2$ agree on their common domain, and $(\forall n \in \omega)(\exists f \in G)(\mathrm{dom}(f)=n+1)$. Therefore, $g= \bigcup G$ is a function with $\mathrm{dom}(g)= \omega$ such that for all $n \in \omega$, $\Psi_{p+1}(n, g \upharpoonright (n+1), V_\alpha)$ holds. This completes the induction step for (II). The base case of the induction on $m$ ($m=1$) follows from the same arguments used to prove the induction step with Lemma \ref{Th:BaseCaseIIKeyLemma} replacing the induction hypothesis. This completes the proof of the lemma.
\Square
\end{proof}

Using Lemma \ref{Th:KeyLemmaForProvingConsistency} we can show that for $m \geq 1$, $\mathbf{M}+ \Pi_{m+1}\textrm{-collection}+\Sigma_{m+2}\textrm{-induction on } \omega$ proves that there exists a transitive model of $\mathrm{Z}+\Pi_m\textrm{-collection}$.

\begin{Theorems1} \label{Th:CollectionPlusInductionProvesConsistency}
Let $m \geq 1$. The theory $\mathbf{M}+\Pi_{m+1}\textrm{-collection}+\Sigma_{m+2}\textrm{-induction on } \omega$ proves that there exists a transitive models of $\mathrm{Z}+\Pi_m\textrm{-collection}$. 
\end{Theorems1}

\begin{proof}
Work in the theory $\mathbf{M}+\Pi_{m+1}\textrm{-collection}+\Sigma_{m+2}\textrm{-induction on } \omega$. By Lemma \ref{Th:KeyLemmaForProvingConsistency}(II), there exists an $f$ such that $\mathrm{dom}(f)= \omega$, and for all $n \in \omega$, $f \upharpoonright (n+1)$ is an $n$-good $\langle m+1, \omega\rangle$-submodel approximation. Let $\beta$ be an ordinal such that $V_\beta= \bigcup \mathrm{rng}(f)$. We claim that $\langle V_\beta, \in \rangle$ is a set structure that satisfies $\mathrm{Z}+\Pi_m\textrm{-collection}$. Since $\beta$ is a limit ordinal $>\omega$, it is immediate that $\langle V_\beta, \in \rangle$ satisfies all of the axioms of $\mathrm{Z}$. Let $\phi(x, y, \vec{z})$ be a $\Pi_m$-formula. Let $\vec{a}, b \in V_\beta$. Note that Definition \ref{Df:nGoodSubmodelApproximation} implies that $V_\beta$ is a $\Pi_{m+1}$-elementary submodel of the universe, and for all $n \in \omega$, $\langle f(n), \in \rangle \prec_m \langle V_\beta, \in \rangle$. Let $k \in \omega$ be such that $\vec{a}, b \in f(k)$. Now, it follows from Definition \ref{Df:nGoodSubmodelApproximation}(v) that for all $x \in b$,
$$\langle V_\beta, \in \rangle \models \exists y \phi(x, y, \vec{a}) \textrm{ if and only if } \langle V_\beta, \in \rangle \models (\exists y \in f(k+1)) \phi^{\langle f(k+1), \in \rangle}(x, y, \vec{a})$$
$$\textrm{if and only if } \langle V_\beta, \in \rangle \models (\exists y \in f(k+1)) \phi(x, y, \vec{a}).$$ 
Therefore
$$\langle V_\beta, \in \rangle \models (\forall x \in b)(\exists y \phi(x, y, \vec{a}) \Rightarrow (\exists y \in f(k+1)) \phi(x, y, \vec{a}))$$
and so $\langle V_\beta, \in \rangle$ satisfies strong $\Pi_m$-collection. Since $\langle V_\beta, \in \rangle$ is a transitive set structure, we can conclude that $\mathbf{M}+\Pi_{m+1}\textrm{-collection}+\Sigma_{m+2}\textrm{-induction on } \omega$ proves that there exists a transitive model of $\mathrm{Z}+\Pi_m\textrm{-collection}$.     
\Square
\end{proof}

We now turn to generalising Theorem \ref{Th:Delta0PCollectionConsistentWithMOST} to show that for all $m \geq 1$, the theories $\mathbf{M}+\textrm{strong } \Pi_m\textrm{-collection}$ and $\mathbf{M}+\Pi_{m+1}\textrm{-collection}$ have the same consistency strength. The key ingredient for this result will be the fact that if $m \geq 1$ and $\mathcal{M}$ is a model of $\mathbf{M}+\textrm{strong } \Pi_m\textrm{-collection}$, then for every standard natural number $n$, there exists an $n$-good $\langle m+1, \omega\rangle$-submodel approximation in $\mathcal{M}$.    

\begin{Lemma1} \label{Th:GeneralFiniteSatisfiabilityLemma}
Let $m \geq 1$ and let $\mathcal{M} \models \mathbf{M}+\textrm{strong } \Pi_m\textrm{-collection}$. For all $n \in \omega$ and for all $\alpha \in \mathrm{Ord}^\mathcal{M}$,
$$\mathcal{M} \models \exists f(f \textrm{ is an }n\textrm{-good } \langle m+1, \alpha \rangle\textrm{-submodel approximation}).$$
\end{Lemma1}

\begin{proof}
Let $\alpha \in \mathrm{Ord}^\mathcal{M}$. We prove the lemma by external induction on $n$. It is clear that 
$$\mathcal{M} \models \exists f(f \textrm{ is a } 0\textrm{-good } \langle m+1, \alpha\rangle\textrm{-submodel approximation}).$$
Suppose that $p \in \omega$ and $f \in M$ are such that 
$$\mathcal{M} \models (f \textrm{ is a } p\textrm{-good } \langle m+1, \alpha\rangle\textrm{-submodel approximation}).$$
Work inside $\mathcal{M}$. Let $V_\delta$ be the rank such that $f(p)=V_\delta$. Consider the $\Pi_m$-formula $\psi(x, y)$ defined by
$$(x=\langle a, l \rangle) \land (l= \ulcorner \phi(u, v) \urcorner \textrm{ where } \phi \textrm{ is } \Pi_m) \land \mathrm{Sat}_{\Pi_m}(l, \langle y, a \rangle).$$
Strong $\Pi_m$-collection implies that there is a set $C$ such that 
$$(\forall x \in V_\delta \times \omega)(\exists y \psi(x, y) \Rightarrow (\exists y \in C)\psi(x, y)).$$
Let $\gamma > \delta$ be such that $C \subseteq V_\gamma$. Using Lemma \ref{Th:BaseCaseIIKeyLemma} (if $m=1$) or Lemma \ref{Th:KeyLemmaForProvingConsistency} (if $m > 1$), we can find a function $g$ with $\mathrm{dom}(g)= \omega$ such that for all $k \in \omega$,
$$g \upharpoonright (k+1) \textrm{ is a } k\textrm{-good } \langle m, \gamma \rangle \textrm{-submodel approximation}.$$
Let $\beta$ be such that $V_\beta= \bigcup \mathrm{rng}(g)$. It follows that for all $l \in \omega$ and for all $a \in V_\beta$,
$$\textrm{if }\langle V_\beta, \in \rangle \models \mathrm{Sat}_{\Pi_m}(l, a), \textrm{ then } \mathrm{Sat}_{\Pi_m}(l, a).$$
And, since $C \subseteq V_\beta$, for all $l \in \omega$ and for all $a \in V_\delta$,
$$\textrm{if } \mathrm{Sat}_{\Sigma_{m+1}}(l, a), \textrm{ then } \langle V_\beta, \in \rangle \models \mathrm{Sat}_{\Sigma_{m+1}}(l, a).$$
Therefore, $h= f \cup \{\langle p+1, V_\beta\rangle\}$ is a $p+1$-good $\langle m+1, \alpha \rangle$-submodel approximation. This concludes the proof of the induction step and the lemma.     
\Square
\end{proof}

We now use a generalisation of the construction used is the proof of Theorem \ref{Th:Delta0PCollectionConsistentWithMOST} to obtain a model $\mathbf{M}+\Pi_{m+1}\textrm{-collection}$ from a model of $\mathbf{M}+\textrm{strong }\Pi_{m}\textrm{-collection}$.

\begin{Theorems1} \label{Th:ConsistencyOfCollectionWithStrongCollection}
Let $m \geq 1$.
\begin{itemize}
\item[(I)] The theory $\mathbf{M}+\Pi_{m+1}\textrm{-collection}$ is $\Pi_{m+3}$-conservative over the theory\\
$\mathbf{M}+\textrm{strong }\Pi_m\textrm{-collection}$. 
\item[(II)] The theory $\mathrm{Mac}+\Pi_{m+1}\textrm{-collection}$ is $\Pi_{m+3}$-conservative over the theory $\mathrm{Mac}+\textrm{strong }\Pi_m\textrm{-collection}$. 
\end{itemize}  
\end{Theorems1}

\begin{proof}
To prove (I) it is sufficient to show that every $\Sigma_{m+3}$-sentence that is consistent with $\mathbf{M}+\textrm{strong }\Pi_m\textrm{-collection}$ is also consistent with $\mathbf{M}+\Pi_{m+1}\textrm{-collection}$. Suppose that $\exists \vec{x} \forall \vec{y} \theta(\vec{x}, \vec{y})$, where $\theta(\vec{x}, \vec{y})$ is a $\Sigma_{m+1}$-formulae, is consistent with $\mathbf{M}+\textrm{strong }\Pi_m\textrm{-collection}$. Let $\mathcal{M}= \langle M, \in^\mathcal{M} \rangle$ be a recursively saturated model of $\mathbf{M}+\textrm{strong }\Pi_m\textrm{-collection}+\exists \vec{x} \forall \vec{y} \theta(\vec{x}, \vec{y})$. Let $\vec{a} \in M$ be such $\mathcal{M}\models \forall \vec{y} \theta(\vec{a}, \vec{y})$ and let $\alpha \in M$ be an ordinal such that $\vec{a} \in (V_\alpha^{\mathcal{M}})^*$. Consider the type
{\small $$\Xi(x, u)= \{x \in \omega\}\cup\{x > n \mid n \in \omega\}\cup \{\exists f(f \textrm{ is an }x\textrm{-good } \langle m+1, \alpha\rangle\textrm{-submodel approximation})\}.$$}
By Lemma \ref{Th:GeneralFiniteSatisfiabilityLemma}, $\Xi(x, u)$ is finitely realised in $\mathcal{M}$, and so there exists $k \in M$ such that $\Xi(k, u)$ is satisfied in $\mathcal{M}$. Note that $k$ is a nonstandard element of $\omega^{\mathcal{M}}$. Let $f \in M$ be such that 
$$\mathcal{M} \models (f \textrm{ is a }k\textrm{-good } \langle m+1, \alpha\rangle\textrm{-submodel approximation}).$$
Define $\mathcal{N}= \langle N, \in^\mathcal{N} \rangle$ by
$$N= \bigcup_{n \in \omega} f(n^\mathcal{M})^* \textrm{ and } \in^\mathcal{N} \textrm{ is the restriction of } \in^\mathcal{M} \textrm{ to }N.$$
We claim that $\mathcal{N}$ satisfies $\mathbf{M}+\Pi_{m+1}\textrm{-collection}+\exists \vec{x} \forall \vec{y} \theta(\vec{x}, \vec{y})$. Note that $\mathcal{N} \subseteq_e^\mathcal{P} \mathcal{M}$. It follows from the fact that $f$ is an $k$-good $\langle m+1, \alpha \rangle$-submodel approximation that $\mathcal{N} \models \mathbf{M}$ and for all $x \in N$, $\mathcal{P}^\mathcal{N}(x)=\mathcal{P}^\mathcal{M}(x)$. Moreover, Definition \ref{Df:nGoodSubmodelApproximation}(iv) implies that $\mathcal{N} \prec_{m+1} \mathcal{M}$. Therefore, since $\vec{a} \in N$, $\mathcal{N}\models \exists \vec{x} \forall \vec{y} \theta(\vec{x}, \vec{y})$.

We are left to show that $\Pi_{m+1}$-collection holds in $\mathcal{N}$. Using exactly the same reasoning that was used in the proof of Theorem \ref{Th:Delta0PCollectionConsistentWithMOST}, we can see that if $C \in M$ is such that $C^* \subseteq N$, then $C \in N$. Now, let $\phi(x, y, \vec{z})$ be a $\Pi_{m+1}$-formula. Let $\vec{d}, b \in N$ be such that
$$\mathcal{N}\models (\forall x \in b) \exists y \phi(x, y, \vec{d})$$
The following formula is a $\mathbf{Bol}(\Pi_{m+1})$-formula with parameters $\vec{d}$, $k$ and $f$:
$$\phi(x, y, \vec{d}) \land (\forall n \in k)(y \notin f(n) \Rightarrow \neg (\exists w \in f(n))\phi(x, w, \vec{d})).$$
And, since $\mathcal{N} \prec_{m+1} \mathcal{M}$,
$$\mathcal{M} \models (\forall x \in b)(\exists y \in f(k))(\phi(x, y, \vec{d}) \land (\forall n \in k)(y \notin f(n) \Rightarrow \neg (\exists w \in f(n))\phi(x, w, \vec{d}))).$$
Working inside $\mathcal{M}$, $\mathbf{Bol}(\Pi_{m+1})$-separation (Lemma \ref{Th:BasicRelationships}) implies that
$$C= \{\langle x, y \rangle \in b \times f(k) \mid \phi(x, y, \vec{d}) \land (\forall n \in k)(y \notin f(n) \Rightarrow \neg (\exists w \in f(n))\phi(x, w, \vec{d}))$$
is a set. And, the fact that $\mathcal{N} \prec_{m+1} \mathcal{M}$ ensures that $C^* \subseteq N$. Therefore $C \in N$. Working inside $\mathcal{N}$, let $B= \mathrm{rng}(C)$. So, 
$$\mathcal{N} \models (\forall x \in b)(\exists y \in B) \phi(x, y, \vec{d}),$$
which shows that $\mathcal{N} \models \Pi_{m+1}\textrm{-collection}$.
 
To see that (II) holds observe that if the Axiom of Choice holds in $\mathcal{M}$ in the proof of (I), then it also holds in $\mathcal{N}$.
\end{proof}

\begin{Coroll1} \label{Th:ConsistencyOfCollectionWithStrongCollection2}
If $\mathbf{M}+\textrm{strong }\Pi_m\textrm{-collection}$ is consistent, then so is $\mathbf{M}+\Pi_{m+1}\textrm{-collection}$.
\end{Coroll1}

Theorem \ref{Th:CollectionPlusInductionProvesConsistency} and Corollary \ref{Th:ConsistencyOfCollectionWithStrongCollection2} yield:

\begin{Coroll1}
If $m \geq 1$, then
$$\mathbf{M}+ \Pi_{m+1}\textrm{-collection} \vdash \mathrm{Con}(\mathbf{M}+\Pi_m\textrm{-collection})$$
\Square
\end{Coroll1}

These results also reveal the limitations of the theory $\mathbf{M}+\Pi_m\textrm{-collection}$ when $m \geq 2$.

\begin{Coroll1} \label{Th:LimitationsOfCollection}
If $m \geq 1$, then
$$\mathbf{M}+\Pi_{m+1}\textrm{-collection} \nvdash \Sigma_{m+2}\textrm{-induction on } \omega.$$
\end{Coroll1}

\begin{proof}
One can easily verify that by starting with a model of $\mathbf{M}+\textrm{strong }\Pi_m\textrm{-collection}+\neg\mathrm{Con}(Z+\Pi_m\textrm{-collection})$ in the proof of Theorem \ref{Th:ConsistencyOfCollectionWithStrongCollection}, one obtains a model of $\mathrm{M}+\Pi_{m+1}\textrm{-collection}+\neg\mathrm{Con}(Z+\Pi_m\textrm{-collection})$. If $\mathbf{M}+\Pi_{m+1}\textrm{-collection}$ proves $\Sigma_{m+2}$-induction, then, by Theorem \ref{Th:CollectionPlusInductionProvesConsistency}, this model would also satisfy $\mathrm{Con}(Z+\Pi_m\textrm{-collection})$, which is a contradiction.
\Square
\end{proof}

The proof of Proposition 9.20 of \cite{mat01} shows that there is an instance of $\Sigma_2$-induction on $\omega$ that coupled with the theory $\mathbf{M}$ proves the consistency of $\mathrm{Mac}$. Therefore, by observing that the proof of Theorem \ref{Th:Delta0PCollectionConsistentWithMOST} can be used to obtain a model of $\mathrm{MOST}+\Pi_1\textrm{-collection}+\neg\mathrm{Con}(\mathrm{MOST})$, we can see that there is an instance of $\Sigma_2$-induction on $\omega$ that is not provable in $\mathrm{MOST}+\Pi_1\textrm{-collection}$. Therefore Corollary \ref{Th:LimitationsOfCollection} also holds when $m=0$.

\section[The strength of $\Pi_n$-collection over $\mathrm{KPI}$]{The strength of $\Pi_n$-collection over $\mathrm{KPI}+V=L$} \label{Sec:ResultsForKP}

In this section we show that the techniques developed in sections \ref{Sec:BaseCaseSection} and \ref{Sec:GeneralCollectionResults} can be adapted to reveal the relative strengths of fragments of the collection scheme over the base theory $\mathrm{KPI}+V=L$. This is achieved by replacing the levels of the $V$-hierarchy in Definition \ref{Df:nGoodSubmodelApproximation} by levels of the $L$-hierarchy. 

\begin{Definitions1} \label{Df:nGoodLApproximation}
Let $n, m \in \omega$, and let $\alpha$ be an ordinal. We say that $f$ is an $n$-good $\langle m+1, \alpha \rangle$-$L$-approximation if 
\begin{itemize}
\item[(i)] $f$ is a function and $\mathrm{dom}(f)= n+1$ 
\item[(ii)] $f(\emptyset)= L_\alpha$
\item[(iii)] $(\forall k \in n+1)\exists \beta((\beta \textrm{ is an ordinal})\land f(k)= L_\beta)$
\item[(iv)]
$$(\forall k \in n)(\forall l \in \omega)(\forall a \in f(k+1))((\langle f(k+1), \in\rangle \models \mathrm{Sat}_{\Pi_m}(l, a)) \Rightarrow \mathrm{Sat}_{\Pi_m}(l, a))$$
\item[(v)]
$$(\forall k \in n)(\forall l \in \omega)(\forall a \in f(k))(\mathrm{Sat}_{\Sigma_{m+1}}(l, a) \Rightarrow (\langle f(k+1), \in  \rangle \models \mathrm{Sat}_{\Sigma_{m+1}}(l, a)))$$
\end{itemize}
\end{Definitions1}

Note that the only difference between Definitions \ref{Df:nGoodSubmodelApproximation} and \ref{Df:nGoodLApproximation} are that the references to levels of the $V$-hierarchy in clauses (ii) and (iii) of Definition \ref{Df:nGoodSubmodelApproximation} have been replaced by level of the $L$-hierarchy in Definition \ref{Df:nGoodLApproximation}. It should be clear that the expression ``$f(\emptyset)= L_\alpha$" remains $\Delta_0$ with parameters $f$ and $L_\alpha$, and, in light of Theorem \ref{Th:DefinabilityOfLInKPI}, the expression ``$(\forall k \in n+1)\exists \beta((\beta \textrm{ is an ordinal})\land f(k)= L_\beta)$" is equivalent to a $\Sigma_1$-formula with parameters $f$ and $n$ in the theory $\mathrm{KPI}$. As we did in section \ref{Sec:GeneralCollectionResults}, we introduce specific notion for formulae that express that $f$ is an $n$-good $\langle m+1, \alpha\rangle$-$L$-approximation.

\begin{Definitions1}
Let $\alpha$ be an ordinal and let $m \in \omega$. We write $\Psi^*_m(n, f, L_\alpha)$ for the formula, with free variables $f$ and $n$, and parameter $L_\alpha$, that the theory $\mathrm{KPI}$ proves asserts that $f$ in an $n$-good $\langle m+1, \alpha\rangle$-$L$-approximation, and such that $\Psi^*_0(n, f, L_\alpha)$ is $\mathbf{Bol}(\Sigma_2)$, and if $m > 0$, $\Psi^*_m(n, f, L_\alpha)$ is $\Pi_{m+1}$. 
\end{Definitions1}

Using the same arguments as we used in the proofs of Lemmas \ref{Th:BaseCaseIKeyLemma} and \ref{Th:BaseCaseIIKeyLemma} we obtain:

\begin{Lemma1} \label{Th:BaseCaseLemmaKPI1}
The theory $\mathrm{KPI}+V=L+\Pi_1\textrm{-collection}+\Sigma_2\textrm{-induction on } \omega$ proves that for all ordinals $\alpha$ and for all $n \in \omega$, there exists an $n$-good $\langle 1, \alpha\rangle$-$L$-approximation.
\Square
\end{Lemma1}

\begin{Lemma1} \label{Th:BaseCaseLemmaKPI2}
The theory $\mathrm{KPI}+V=L+\Pi_1\textrm{-collection}+\Sigma_2\textrm{-induction on } \omega$ proves that for all ordianls $\alpha$, there exists a function $f$ with $\mathrm{dom}(f)=\omega$ such that for all $n \in \omega$, $f \upharpoonright (n+1)$ in an $n$-good $\langle 1, \alpha\rangle$-$L$-approximation.
\Square
\end{Lemma1}

Lemmas \ref{Th:BaseCaseLemmaKPI1} and \ref{Th:BaseCaseLemmaKPI2} now provide the base case of an induction argument that proves an analogue of Lemma \ref{Th:KeyLemmaForProvingConsistency}.

\begin{Lemma1} \label{Th:KeyLemmaForProvingConsistencyKPI}
Let $m \in \omega$. The theory $\mathrm{KPI}+V=L+\Pi_{m+1}\textrm{-collection}+\Sigma_{m+2}\textrm{-induction}$ proves
\begin{itemize}
\item[(I)] for all ordinals $\alpha$ and for all $n \in \omega$, there exists an $n$-good $\langle m+1, \alpha\rangle$-$L$-approximation,
\item[(II)] for all ordinals $\alpha$, there exists a function $f$ with $\mathrm{dom}(f)= \omega$ such that for all $n \in \omega$, $f \upharpoonright (n+1)$ is an $n$-good $\langle m+1, \alpha\rangle$-$L$-approximation.  
\end{itemize}
\Square
\end{Lemma1}

Lemma \ref{Th:KeyLemmaForProvingConsistencyKPI} provides the key ingredient for showing that the theory\\ $\mathrm{KPI}+V=L+\Pi_{m+1}\textrm{-collection}+\Sigma_{m+2}\textrm{-induction on } \omega$ proves the consistency of the theory $\mathrm{KPI}+V=L+\textrm{strong } \Pi_m\textrm{-collection}+\Pi_\infty\textrm{-foundation}$.

\begin{Theorems1} \label{Th:ProofOfCosistencyForKPI}
Let $m \in \omega$. The theory $\mathrm{KPI}+V=L+\Pi_{m+1}\textrm{-collection}+\Sigma_{m+2}\textrm{-induction on } \omega$ proves that there exists a transitive model of $\mathrm{KPI}+V=L+\textrm{strong }\Pi_m\textrm{-collection}+\Pi_\infty\textrm{-foundation}$.
\end{Theorems1}

\begin{proof}
Work in the theory $\mathrm{KPI}+V=L+\Pi_{m+1}\textrm{-collection}+\Sigma_{m+2}\textrm{-induction on } \omega$. By Lemma \ref{Th:KeyLemmaForProvingConsistencyKPI}(II), there exists $f$ such that $\mathrm{dom}(f)=\omega$, and for all $n \in \omega$, $f\upharpoonright (n+1)$ is an $n$-good $\langle m+1, \omega\rangle$-$L$-approximation. Let $\beta$ be an ordinal such that $L_\beta= \bigcup \mathrm{rng}(f)$. We claim that $\langle L_\beta, \in \rangle$ is a set structure that satisfies $\mathrm{KPI}+\textrm{strong }\Pi_m\textrm{-collection}+\Pi_\infty\textrm{-foundation}$ (=$\mathbf{M}^-+\textrm{strong }\Pi_m\textrm{-collection}+\Pi_\infty\textrm{-foundation}$). Note that, since $\beta$ is a limit ordinal, $L_\beta$ is a transitive set that is closed under G\"{o}del operations. Therefore $\langle L_\beta, \in\rangle$ satisfies all of the axioms of $\mathbf{M}^-$. Let $\phi(x, \vec{z})$ be a $\Pi_\infty$-formula and let $\vec{a} \in L_\beta$. Separation in the theory $\mathrm{KPI}$ implies that 
$$A= \{ x \in L_\beta \mid \langle L_\beta, \in \rangle \models \phi(x, \vec{a})\}$$
is a set. Therefore, set foundation in $\mathrm{KPI}$, implies that if $A\neq \emptyset$, then $A$ has an $\in$-least element. This shows that $\langle L_\beta, \in\rangle$ satisfies $\Pi_\infty$-foundation. Finally, identical reasoning to that used in the proof of Theorem \ref{Th:CollectionPlusInductionProvesConsistency} shows that $\langle L_\beta, \in\rangle$ satisfies $\textrm{strong }\Pi_m\textrm{-collection}$. Since $\langle L_\beta, \in\rangle$ is a transitive set structure, we can conclude that $\mathrm{KPI}+\Pi_{m+1}\textrm{-collection}+\Sigma_{m+2}\textrm{-induction on } \omega$ proves that there exists a transitive models of $\mathrm{KPI}+\textrm{strong }\Pi_m\textrm{-collection}+\Pi_\infty\textrm{-foundation}+V=L$. 
\Square
\end{proof}

We next turn indicating how the proof of Theorem \ref{Th:ConsistencyOfCollectionWithStrongCollection} can be adapted to obtain an analogue of this result with the base theory $\mathbf{M}$ replaced by $\mathrm{KPI}+V=L$. The same argument used in the proof of Lemma \ref{Th:GeneralFiniteSatisfiabilityLemma} can be used to prove the following:

\begin{Lemma1} \label{Th:FiniteSatisfiabilityLemmaKPI}
Let $m \in \omega$ and let $\mathcal{M} \models\mathrm{KPI}+V=L+\textrm{strong }\Pi_m\textrm{-collection}$. For all $n \in \omega$ and for all $\alpha \in \mathrm{Ord}^\mathcal{M}$,
$$\mathcal{M} \models \exists f(f \textrm{ is an } n \textrm{-good } \langle m+1, \alpha \rangle\textrm{-}L\textrm{-approximation}).$$
\Square
\end{Lemma1}

Lemma \ref{Th:FiniteSatisfiabilityLemmaKPI} yields an analogue of Theorem \ref{Th:ConsistencyOfCollectionWithStrongCollection}.

\begin{Theorems1}\label{Th:ConsistencyOfStrongCollectionWithCollectionKPI}
Let $m \in \omega$. 
\begin{itemize}
\item[(I)] The theory $\mathrm{KPI}+V=L+\Pi_{m+1}\textrm{-collection}$ is $\Pi_{m+3}$-conservative over the theory $\mathrm{KPI}+V=L+\textrm{strong }\Pi_m\textrm{-collection}$. 
\item[(II)] If $\mathrm{KPI}+V=L+\textrm{strong }\Pi_m\textrm{-collection}$ is consistent, then so is $\mathrm{KPI}+V=L+\Pi_{m+1}\textrm{-collection}$.
\end{itemize}
\Square 
\end{Theorems1}

Theorems \ref{Th:ProofOfCosistencyForKPI} and \ref{Th:ConsistencyOfStrongCollectionWithCollectionKPI} yield:

\begin{Coroll1}
If $m \geq 1$, then
$$\mathrm{KPI}+V=L+\Pi_{m+1}\textrm{-collection} \vdash \mathrm{Con}(\mathrm{KPI}+V=L+\Pi_m\textrm{-collection})$$
\Square
\end{Coroll1}

\begin{Quest1}
Does the theory $\mathrm{KPI}+V=L+\textrm{strong }\Pi_0\textrm{-collection}$ prove the consistency of $\mathrm{KPI}$?
\end{Quest1}

I am grateful to Ali Enayat for the following observation:

\begin{Remark1}\label{Th:AliRemark}
The proofs of Theorems \ref{Th:Delta0PCollectionConsistentWithMOST}, \ref{Th:Delta0PCollectionConsistentWithMac}, \ref{Th:ConsistencyOfCollectionWithStrongCollection} and \ref{Th:ConsistencyOfStrongCollectionWithCollectionKPI} can all be formalised in the subsystem of second order arithmetic $\mathrm{WKL}_0$. The fact that $\mathrm{WKL}_0$ is conservative over Primitive Recursive Arithmetic ($\mathrm{PRA}$) for sentences that are $\Pi_2$ sentences of arithmetic (see \cite[Theorem IX.3.16]{sim09}), then shows that all of these results are theorems of $\mathrm{PRA}$.  
\end{Remark1}

\noindent{\bf Acknowledgements:} I am very grateful to Adrian Mathias and Ali Enayat for their helpful comments on earlier drafts of this paper. In particular, Ali Enayat's observations led to the strengthening of Theorems \ref{Th:Delta0PCollectionConsistentWithMOST}, \ref{Th:Delta0PCollectionConsistentWithMac}, \ref{Th:ConsistencyOfCollectionWithStrongCollection} and \ref{Th:ConsistencyOfStrongCollectionWithCollectionKPI}. I would also like to thank the anonymous referee for their careful reading of this paper and their thoughtful suggestions.

\bibliographystyle{alpha}
\bibliography{.}          

\begin{thebibliography}{9}

\bibitem[Bar]{bar75} Barwise, Jon. \emph{Admissible Sets and Structures}. Perspectives in Mathematical Logic. Springer-Verlag. 1975.

\bibitem[FK]{fk91} Forster, Thomas and Kaye, Richard. ``End-extensions preserving power set". \emph{The Journal of Symbolic Logic}. Vol. 56. No. 1. 1991. pp 323--328.

\bibitem[Fri]{fri73} Friedman, Harvey M. ``Countable models of set theories". \emph{Cambridge Summer School in Mathematical Logic, August 1--21, 1971}. Edited by A. R. D. Mathias and H. Rogers Jr. Springer Lecture Notes in Mathematics. Vol. 337. Springer, Berlin. 1973. pp 539--573.

\bibitem[FLW]{flw16} Friedman, Sy-David; Li, Wei and Wong, Tin Lok. ``Fragments of Kripke-Platek Set Theory and the Metamathematics of $\alpha$-Recursion Theory". \emph{Archive for Mathematical Logic}. Vol. 55. No. 7. 2016. pp 899--924.  

\bibitem[M]{mck15} McKenzie, Zachiri. ``Automorphisms of models of set theory and extensions of $\mathrm{NFU}$". \emph{Annals of Pure and Applied Logic}. Vol. 166. 2015. pp 601--638.

\bibitem[Mat69]{mat69} Mathias, Adrian R. D. ``Notes on set theory". Available online: \texttt{https://www.dpmms.cam.ac.uk/\textasciitilde ardm/} (last accessed on 29/vi/2018)

\bibitem[Mat01]{mat01} Mathias, Adrian R. D. ``The strength of Mac Lane set theory". \emph{Annals of Pure and Applied Logic}. Vol. 110. 2001. pp 107--234.

\bibitem[Res]{res87} Ressayre, Jean-Pierre. ``Mod\`{e}les non standard et sous-syst\`{e}mes remarquables de ZF". \emph{Mod\`{e}les non standard en arithm\'{e}tique et th\'{e}orie des ensembles}. Publications
Math\'{e}matiques de l'Universit\'{e} Paris VII. Vol. 22. Universit\'{e} de Paris VII, U.E.R. de Math\'{e}matiques, Paris. 1987. pp 47--147.

\bibitem[Sim]{sim09} Simpson, Stephen G. \emph{Subsystems of second order arithmetic}. Association for Symbolic Logic Perspectives in Logic series. Cambridge University Press. 2009.

\bibitem[Tak]{tak72} Takahashi, Moto-o. ``$\tilde{\Delta}_1$-definability in set theory". \emph{Conference in mathematical logic --- London '70}. Edited by W. Hodges. Springer Lecture Notes in Mathematics. Vol. 255. Springer. 1972. pp 281--304.

\end{thebibliography}
    
\end{document}